\documentclass[10pt]{article}

\usepackage{amsmath}
\usepackage{amsfonts}
\usepackage{amssymb}
\usepackage{framed}
\usepackage{enumerate}
\usepackage{mathrsfs}
\usepackage[hidelinks]{hyperref}
\usepackage{enumitem}
\usepackage{mathrsfs}
\usepackage{scrextend}
\usepackage{amsthm}
\usepackage{bbm}
\usepackage{thmtools}
\usepackage{thm-restate}
\usepackage{mathabx}
\usepackage[margin=1.5in]{geometry}
\usepackage{rotating}

\usepackage{xcolor}
\hypersetup{
    colorlinks,
    linkcolor={red!50!black},
    citecolor={blue!50!black},
    urlcolor={blue!80!black}
}

\usepackage{tikz-cd}

\theoremstyle{plain}
\newtheorem{thm}{Theorem}[section]
\newtheorem{lem}[thm]{Lemma}
\newtheorem{prop}[thm]{Proposition}
\newtheorem{cor}[thm]{Corollary}

\theoremstyle{definition}
\newtheorem{defn}[thm]{Definition}

\newtheorem{claim}[thm]{Claim}

\theoremstyle{remark}
\newtheorem*{rem}{Remark}

\newtheorem*{notn}{Notation}

\tikzset{
  symbol/.style={
    draw=none,
    every to/.append style={
      edge node={node [sloped, allow upside down, auto=false]{$#1$}}}
  }
}

\newcommand{\Spec}{\textrm{Spec} \hspace{0.15em} }

\newcommand\restr[2]{{
	\left.\kern-\nulldelimiterspace
	#1
	\vphantom{\big|}
	\right|_{#2}
	}}
\newcommand{\an}[1]{#1^{\textrm{an}}}

\newcommand{\FSch}{\textrm{FSch}}
\newcommand{\Set}{\textrm{Set}}
\newcommand{\Sch}{\textrm{Sch}}
\newcommand{\An}{\textrm{An}}

\newcommand{\Hom}{\textrm{Hom}}

\newcommand{\ch}[1]{\widecheck{{#1}}}

\newcommand{\codim}{\textrm{codim}}

\newcommand{\GL}{\textrm{GL}}

\newcommand{\Qbar}{\overline{\mathbb{Q}}}

\usepackage{amsmath,calligra,mathrsfs}

\DeclareMathOperator{\sheafhom}{\mathscr{H}\text{\kern -3pt {\calligra\large om}}\,}

\title{On the Transcendence of Period Images}
\author{David Urbanik}

\begin{document}

\maketitle

\begin{abstract}
Let $f : X \to S$ be a family of smooth projective algebraic varieties over a smooth connected base $S$, with everything defined over $\Qbar$. Denote by $\mathbb{V} = R^{2i} f_{*} \mathbb{Z}(i)$ the associated integral variation of Hodge structure on the degree $2i$ cohomology. We consider the following question: when can a fibre $\mathbb{V}_{s}$ above an algebraic point $s \in S(\Qbar)$ be isomorphic to a transcendental fibre $\mathbb{V}_{s'}$ with $s' \in S(\mathbb{C}) \setminus S(\Qbar)$? When $\mathbb{V}$ induces a quasi-finite period map $\varphi : S \to \Gamma \backslash D$, conjectures in Hodge theory predict that such isomorphisms cannot exist. We introduce new differential-algebraic techniques to show this is true for all points $s \in S(\Qbar)$ outside of an explicit proper closed algebraic subset of $S$. As a corollary we establish the existence of a canonical $\overline{\mathbb{Q}}$-algebraic model for normalizations of period images.
\end{abstract}

\tableofcontents

\section{Introduction}

Suppose that $f : X \to S$ is a family of smooth projective algebraic varieties over a smooth connected base $S$, all defined over a number field $K \subset \mathbb{C}$. Fix a cohomological degree $2i$, and let $\mathbb{V} = R^{2i} \an{f}_{*} \mathbb{Z}(i)$ be the variation of Hodge structure determined by $f$. It is in general a difficult question to relate the arithmetic structure of the family $f$ to the transcendental object $\mathbb{V}$, and many central conjectures in Hodge theory are devoted to elucidating this relationship. One simple question which motivates our paper is: when can the Hodge structure $\mathbb{V}_{s}$ with $s \in S(\Qbar)$ be isomorphic to a Hodge structure $\mathbb{V}_{s'}$ with $s' \in S(\mathbb{C}) \setminus S(\Qbar)$?

Let us first explain the conjectural solution to this problem. The Hodge conjecture (or weaker versions, such as absolute Hodge) predicts that if $\mathbb{V}_{s}$ and $\mathbb{V}_{s'}$ are isomorphic, then so are $\mathbb{V}_{s^{\sigma}}$ and $\mathbb{V}_{s'^{\sigma}}$ for any $\sigma \in \textrm{Aut}(\mathbb{C}/K)$. This shows that the equivalence relation (E1) on $S(\mathbb{C})$ defined by $s \sim_{(E1)} s'$ if $\mathbb{V}_{s} \simeq \mathbb{V}_{s'}$ is preserved by conjugation by $\textrm{Aut}(\mathbb{C}/K)$. If the equivalence classes are finite (i.e., the period map $\varphi : S \to \Gamma \backslash D$ induced by $\mathbb{V}$ is quasi-finite) then $s \sim_{(E1)} s'$ implies the orbits of $s$ and $s'$ under $\textrm{Aut}(\mathbb{C}/K)$ have the same cardinality, and so $\sim_{(E1)}$ does not identify algebraic and transcendental points. In the general case where $\varphi$ isn't quasi-finite, we may first complete $S$ analogously to \cite[Cor.~13.7.6]{CMS} to reduce to the case where the map $\varphi$ is proper. Using the Stein factorization of $\varphi$ we then obtain a decomposition $\varphi = \varphi_{\textrm{qf}} \circ \varphi_{\textrm{con}}$, where $\varphi_{\textrm{con}}$ has connected fibres and $\varphi_{\textrm{qf}}$ is quasi-finite, so the problem reduces to understanding the quasi-finite period map $\varphi_{\textrm{qf}}$. For simplicity we assume that $\varphi$ is a proper map for the remainder of the introduction, as in each of our theorems the non-proper case will follow formally from the proper one.

Understanding the arithmetic behaviour of the equivalence relation (E1) seems to be a difficult problem in general. Our central observation in this work is that we can make considerable progress by replacing (E1) with two closely related equivalence relations (E2) and (E3) which are much easier to understand arithmetically, and use this to recover arithmetic information about (E1). To do this, we recall the typical setup for studying period maps. We fix an integral polarized lattice $(V, Q)$ of the same type as the variation $\mathbb{V}$, denote by $\ch{D}$ the flag variety parametrizing Hodge flags on $V$ satisfying the first Hodge-Riemann bilinear relation, and denote by $D \subset \ch{D}$ its analytic open subset of polarized integral Hodge structures. If $G = \textrm{Aut}(V, Q)$ denotes the symmetry group of linear maps preserving the polarization, then $\ch{D}$ is a homogeneous space for $G(\mathbb{C})$. Letting $\Gamma = G(\mathbb{Z})$ the analytic space $\Gamma \backslash D$ may be interpreted as a moduli space for polarized integral Hodge structures on $V$, and we obtain a canonical analytic map $\varphi : S \to \Gamma \backslash D$ by sending a point $s$ to the isomorphism class of $\mathbb{V}_{s}$. We also fix a lift $\widetilde{\varphi} : \widetilde{S} \to D$, with $\pi : \widetilde{S} \to S$ the universal covering, and denote by $\widetilde{T} \subset D$ the image. Using the fact that $\varphi$ is proper, we show in \autoref{closedT} this is a closed analytic subset of $D$.

We now introduce (E2) and (E3). These equivalence relations will relate points $s$ and $s'$ in $S(\mathbb{C})$ in terms of certain local lifts of $\varphi$ defined on balls centred around $s$ and $s'$.

\begin{notn}
For an analytic variety $X$ and a point $x \in X$ we denote by $(X, x)$ the germ of $X$ at $x$.
\end{notn}

\begin{defn}
Given a period map $\varphi : S \to \Gamma \backslash D$ and a point $s \in S(\mathbb{C})$, a \emph{distinguished local lift} at $s$ of $\varphi$ is a map $\psi : B \to D$ with $B \subset S$ an open ball containing $s$ such that:
\begin{itemize}
\item[(i)] $\psi$ is a local lift, i.e., $\pi \circ \psi = \restr{\varphi}{B}$, where $\pi : D \to \Gamma \backslash D$ is the natural projection;
\item[(ii)] if $t = \psi(s)$, then $\psi(B)$ lies in a unique analytic component $\mathcal{C}(\psi, s)$ of the germ $(\widetilde{T}, t)$.
\end{itemize}
\end{defn}

\noindent Given a local lift $\psi : B \to D$ containing $s \in S(\mathbb{C})$, one can always obtain a distinguished local lift by shrinking $B$. We then have: 

\begin{defn}
\label{equivdef}
We say:
\begin{itemize}[leftmargin=.5in]
\item[(E1)] that $s_{1} \sim_{(E1)} s_{2}$ if $\varphi(s_{1}) = \varphi(s_{2})$;
\item[(E2)] that $s_{1} \sim_{(E2)} s_{2}$ if there exists distinguished local lifts $\psi_{i} : B_{i} \to D$ with $s_{i} \in B_{i}$ and $\gamma \in G(\mathbb{Z})$ such that $\gamma \cdot \mathcal{C}(\psi_{1}, s_{1}) = \mathcal{C}(\psi_{2}, s_{2})$;
\item[(E3)] that $s_{1} \sim_{(E3)} s_{2}$ if there exists distinguished local lifts $\psi_{i} : B_{i} \to D$ with $s_{i} \in B_{i}$ and $g \in G(\mathbb{C})$ such that $g \cdot \mathcal{C}(\psi_{1}, s_{1}) = \mathcal{C}(\psi_{2}, s_{2})$.
\end{itemize}
\end{defn}

\vspace{0.5em}
Let us explain the significance of (E2) and (E3), denoting by $T = \varphi(S)$ the image of $\varphi$ in what follows; we note that $T$ is a quasi-projective complex algebraic variety by \cite{OMINGAGA}. The relation (E2) is finer than (E1), but not by much: whereas the equivalence classes of (E1) correspond to the fibres of $\varphi$, we explain in \autoref{E2isnormalization} that the equivalence classes of (E2) are the fibres of the first map in the factorization $S \to T^{\nu} \to T$, where $T^{\nu} \to T$ is the normalization. The normalization map $T^{\nu} \to T$ is bijective exactly over the unibranch locus of $T$, hence (E1) and (E2) agree above the unibranch locus of the period image; in particular, if $T$ is unibranch, then (E1) and (E2) are the same.

The idea is that while replacing (E1) with (E2) does not change the problem much, the derivatives of the distinguished local lifts $\psi_{i} : B_{i} \to D$ which determine the components $\mathcal{C}(\psi_{i}, s_{i})$ for $i = 1, 2$ have arithmetic properties which make (E2) amenable to study by algebraic methods.\footnote{This is roughly analogous to how generic Torelli-type theorems are often proven by replacing the Hodge structure at a point with the infinitesimal variation in a neighbourhood.} In particular it allows us to replace (E2) with the coarser equivalence relation (E3). We will show that for variations of Hodge structure arising from geometry, (E3) can be studied entirely using algebraic methods, leading to the following results:

\begin{defn}
\label{foliationdef}
Given a variety $S$, we say that $S$ admits an algebraic foliation by a set $\mathcal{F}$ of subvarieties of $S$ if there exists a dominant algebraic map $g : S \to U$ such that the (geometrically) irreducible components of the non-empty fibres of $g$ are the elements of $\mathcal{F}$.
\end{defn}

\begin{rem}
Note that in \autoref{foliationdef} we do not assume that the varieties in $\mathcal{F}$ give a partition of $S$, just that the fibres of $g$ do. In other words, the components of the fibres of $g$ are allowed to intersect.
\end{rem}

\begin{defn}
If $Y \subset S$ is a subvariety and $\mathbb{V}$ is a variation of Hodge structure on $S$, we say that $Y$ is \emph{Hermitian symmetric} if the variation $\restr{\mathbb{V}}{Y}$ on $Y$ is Hermitian symmetric (see \autoref{hermdef}).
\end{defn}

\begin{prop}
\label{E3prop}
Let $\varphi : S \to \Gamma \backslash D$ be the period map arising from a $\Qbar$-algebraically-defined variation\footnote{See \autoref{Kalgdef} below; for instance, $\mathbb{V}$ can be induced by a smooth projective $\Qbar$-algebraic family $f : X \to S$.} with $S$ smooth, with $\Gamma = \textrm{Aut}(V,Q)(\mathbb{Z})$, as above. Then
\begin{itemize}
\item[(i)] the graph $R_{3} \subset S \times S$ of the equivalence relation (E3) is a $\Qbar$-algebraic subvariety;
\item[(ii)] either the map
\[ S(\mathbb{C}) \big/ \sim_{(E2)} \hspace{0.5em} \to  \hspace{0.5em} S(\mathbb{C}) \big/ \sim_{(E3)} \]
induced by the identity has finite fibres, or there exists an algebraic foliation $g : S \to U$ by Hermitian symmetric subvarieties $Y \subset S$ whose images $Z = \varphi(Y) \subset T$ are positive-dimensional.
\end{itemize}
\end{prop}

\medskip

Using \autoref{E3prop} (or rather, a slight strengthening of it: we will have a concrete description of the foliation $g$, see \autoref{hermfolprop}), we will then obtain:

\begin{thm}
\label{mainthm}
Let $\mathbb{V}$ be a $\Qbar$-algebraically-defined variation on the smooth $\Qbar$-variety $S$. Let $\varphi : S \to \Gamma \backslash D$ be the associated period map with $\Gamma = \textrm{Aut}(V, Q)(\mathbb{Z})$, which we may assume to be proper after possibly completing $S$. Let $T = \varphi(S)$ be the period image. Then in the factorization $S \xrightarrow{p} T^{\nu} \xrightarrow{\textrm{nor}} T$, where $T^{\nu} \xrightarrow{\textrm{nor}} T$ is the normalization, both $T^{\nu}$ and the map $p$ admit $\overline{\mathbb{Q}}$-algebraic models.
\end{thm}

\begin{cor}
\label{maincor}
In the situation of \autoref{mainthm}, suppose that $\varphi$ is quasi-finite, and let $E \subset S$ be the proper algebraic subvariety of points lying above the non-unibranch locus of $T = \varphi(S)$. Then for any $s \in S(\Qbar) \setminus (E(\mathbb{C}) \cap S(\Qbar))$, there does not exist $s' \in S(\mathbb{C}) \setminus S(\Qbar)$ such that $\mathbb{V}_{s} \simeq \mathbb{V}_{s'}$ as polarized integral Hodge structures.
\end{cor}

\noindent Both \autoref{mainthm} and \autoref{maincor} are consequences of \autoref{E3prop} (more precisely, of \autoref{R3alg} and \autoref{hermfolprop}), as we explain in \autoref{mainresultssec}. Roughly speaking, whenever $S$ does not admit a positive-dimensional Hermitian symmetric foliation we will find that the graph $R_{2} \subset R_{3}$ of (E2) is a union of components of $R_{3}$, hence is defined over $\Qbar$ using \autoref{E3prop}(i). In the general case the Hermitian symmetric part of the variation can be factored out and handled separately, allowing us to argue inductively. This will enable us to put a $\Qbar$-structure on $T^{\nu}$ and $p$, and prove both \autoref{mainthm} and \autoref{maincor}.

In fact, a stronger version of \autoref{mainthm} which replaces $\Gamma = \textrm{Aut}(V, Q)(\mathbb{Z})$ by an arbitrary sublattice $\Gamma' \subset \Gamma$ containing the image of monodromy follows formally from \autoref{mainthm}, see \autoref{anylattice}.

\subsection{Jets and Local Period Maps}

From the preceding discussion, we see that the central focus of the paper is to understand the equivalence relation (E3), and in particular, to establish the assertions (i) and (ii) made in \autoref{E3prop}. A crucial tool will be a new Hodge-theoretic construction associated to a $K$\emph{-algebraically-defined} polarized variation of Hodge structure $\mathbb{V}$, where $K \subset \mathbb{C}$ is a subfield, and the higher dimensional jets determined by \emph{local period maps} associated to $\mathbb{V}$. We now formally introduce both notions, the second of which is implicit in \autoref{equivdef}.

\begin{defn}
\label{Kalgdef}
Let $\mathbb{V}$ be a polarized variation of Hodge structure on the $K$-algebraic variety $S$, with $K \subset \mathbb{C}$. We say that $\mathbb{V}$ is $K$-algebraically-defined if the Hodge bundle $\an{\mathcal{H}} = \mathbb{V} \otimes \mathcal{O}_{\an{S}}$, the Hodge filtration $F^{\bullet} \an{\mathcal{H}}$, the connection $\nabla : \an{\mathcal{H}} \to \an{\mathcal{H}} \otimes \Omega^{1}_{\an{S}}$, and the polarization $\mathcal{Q} : \an{\mathcal{H}} \otimes \an{\mathcal{H}} \to \mathcal{O}_{\an{S}}$ all admit $K$-algebraic models, and these models are compatible with all natural identifications.
\end{defn}

\begin{rem}
As explained in \cite[\S 2]{absladic}, any ``motivic'' variation of Hodge structure, and in particular those arising from geometric $K$-algebraic families $f: X \to S$, is $K$-algebraically-defined. In the case of the variation on primitive cohomology associated to the family $f$, the $K$-algebraicity of the Hodge bundle and filtration comes from the the Hodge-to-deRham spectral sequence, the connection is $K$-algebraic by a result of Katz and Oda \cite{katz1968}, and the polarization is $K$-algebraic due to the compatibility of the algebraic and topological Chern classes and cup products.
\end{rem}

\begin{defn}
Let $\ch{D}$ be the flag variety parametrizing polarized Hodge flags on the lattice $V$ with the same Hodge numbers as the fibres of $\mathbb{V}$, and let $G = \textrm{Aut}(V, Q)$ be the natural symmetry group. A \emph{local period map} associated to $\mathbb{V}$ is then a map $\psi : B \to \ch{D}$, where $B \subset S$ is an analytic ball, which is of the form $\psi = g \cdot  \psi'$, where $g \in G(\mathbb{C})$ and $\psi' : B \to D \subset \ch{D}$ is a local lift of the period map $\varphi : S \to \Gamma \backslash D$ associated to $\mathbb{V}$.
\end{defn}

We will see in \autoref{locpermapredef} that local period maps are maps constructed in the same way as local lifts of $\varphi$ except that instead of identifying the fibres of $\mathbb{V}$ above $B$ using an integral basis of $\mathbb{V}(B)$, we allow ourselves to use more general complex bases in $\mathbb{V}_{\mathbb{C}}(B)$. In particular, there will be a universal $K$-algebraic system of differential equations determined by the connection $\nabla$, the filtration $F^{\bullet}$, and the polarization $\mathcal{Q}$, and the local period maps on $B$ will be exactly the solutions on $B$ to this system. Our main results will require an analysis of the arithmetic properties of this differential system and its solutions, for which we use \emph{jets}.

\begin{notn}
Let $d, r \geq 0$ be integers. We let $\mathbb{D}^{d}_{r} = \Spec K[t_{1}, \hdots, t_{d}]/(t_{1}, \hdots, t_{d})^{r+1}$ be the formal $d$-dimensional disk of order $r$.
\end{notn}

\begin{defn}
Let $S$ be a finite-type algebraic $K$-scheme. For integers $d, r \geq 0$, the \emph{jet space} $J^{d}_{r} S$ is defined to be the $K$-scheme which represents the contravariant functor $\Sch/K \to \Set$ given by 
\[ T \mapsto \Hom_{K}(T \times_{K} \mathbb{D}^{d}_{r}, S), \hspace{1.5em} [T \to T'] \mapsto [\Hom_{K}(T' \times_{K} \mathbb{D}^{d}_{r}, S) \to \Hom_{K}(T \times_{K} \mathbb{D}^{d}_{r}, S)] , \]
where the natural map $\Hom_{K}(T' \times_{K} \mathbb{D}^{d}_{r}, S) \to \Hom_{K}(T \times_{K} \mathbb{D}^{d}_{r}, S)$ obtained by pulling back along $T \times_{K} \mathbb{D}^{d}_{r} \to T' \times_{K} \mathbb{D}^{d}_{r}$. 
\end{defn}

\begin{rem}
In most cases, jets spaces are only studied with $d = 1$. We will find the additional flexibility that comes with using higher-dimensional jets useful, and review the necessary background in \autoref{jetsec}.
\end{rem}

\begin{rem}
We also have an analytic version of this construction, for which we refer to \autoref{anjetspace}.
\end{rem}

\begin{notn}
Given a point $s \in S(\mathbb{C})$, we denote by $(J^{d}_{r} S)_{s}$ the fibre of the natural projection $J^{d}_{r} S \to S$ above $s$.
\end{notn}

\medskip

We will often view points $j$ of $J^{d}_{r} S$ as maps $\mathbb{D}^{d}_{r} \to S$ using representability, and hence given a local period map $\psi : B \to \ch{D}$ we will write $\psi \circ j$ for the image of $j$ in $J^{d}_{r} \ch{D}$ under the map $J^{d}_{r} \psi : J^{d}_{r} B \to J^{d}_{r} \ch{D}$ induced by functoriality. Moreover, the action of $G$ on the flag variety $\ch{D}$ induces an action on $J^{d}_{r} \ch{D}$ for any $d, r \geq 0$, and we will be interested in the quotient $G \backslash J^{d}_{r} \ch{D}$. This is in general only an Artin stack, but this will be enough algebraic structure for our purposes. Our key theoretical result is then the following:

\begin{restatable}{thm}{jetcorresp}
\label{jetcorresp}
Given a $K$-algebraically-defined polarized variation of Hodge structure $\mathbb{V}$ on a smooth $K$-variety $S$, and integers $d, r \geq 0$, there exists a canonical $K$-algebraic (in particular, $\textrm{Aut}(\mathbb{C}/K)$-invariant) map
\[ \eta^{d}_{r} : J^{d}_{r} S \to G \backslash J^{d}_{r} \ch{D} \]
characterized by the following property: for any $j \in (J^{d}_{r} S)(\mathbb{C})$ and any local period map $\psi : B \to \ch{D}$ defined at the base-point of $j$ in $S$, we have $\eta^{d}_{r}(j) = \psi \circ j$ modulo $G$. 
\end{restatable}

The proof of \autoref{E3prop}(i) will actually use a slightly stronger version of \autoref{jetcorresp} which gives a precise description of the $G$-torsor defining the map $\eta^{d}_{r}$. The essential idea is that in most situations, the condition that there exists $g \in G(\mathbb{C})$ such that $g \cdot \mathcal{C}(\psi_{1}, s_{1}) = \mathcal{C}(\psi_{2}, s_{2})$ can be rephrased as a condition on the maps $\psi_{i}$, and this condition will amount to the existence of sufficiently high-order jets $j_{1}$ and $j_{2}$ with $j_{i} \in (J^{d}_{r} S)_{s_{i}}$ such that $\eta^{d}_{r}(j_{1}) = \eta^{d}_{r}(j_{2})$. In particular, the correspondence \autoref{jetcorresp} will allow us to exhibit a countable intersection of $K$-constructible sets $R_{4} \subset R_{3}$ such that $R_{3}$ coincides with its closure.

\subsection{Foliations from Local Symmetries}
\label{folfromloc}

After disposing of \autoref{E3prop}(i), it will remain to deal with \autoref{E3prop}(ii). To illustrate the idea, let us fix a lift $\widetilde{\varphi} : \widetilde{S} \to D \subset \ch{D}$ of the map $\varphi$, and write $\widetilde{T} = \widetilde{\varphi}(\widetilde{S})$ for the image. Let $s \in S(\mathbb{C})$ be a point, and let $t \in \widetilde{T}$ be a point lying above $\varphi(s)$. We consider the subset $H_{t} \subset G(\mathbb{C})$ of elements $g$ mapping $(\widetilde{T}, t)$ to another germ of $\widetilde{T}$; the set $H_{t}$ consists of elements $g$ responsible for identifying $s$ with other points of $S(\mathbb{C})$ under (E3). 

One can check that $H_{t}$ is an analytic set in a neighbourhood of the identity in $G(\mathbb{C})$. In fact, we will prove a stronger result: the germ $(H_{t}, \textrm{id})$ agrees with the germ of a complex Lie subgroup $H \subset G(\mathbb{C})$, and $H$ is independent of $t$. If the map 
\[ S(\mathbb{C}) \big/ \sim_{(E2)} \hspace{0.5em} \to  \hspace{0.5em} S(\mathbb{C}) \big/ \sim_{(E3)} \]
is not quasi-finite, we will argue that $H$ must be positive-dimensional, from which one learns that $\widetilde{T}$ admits a non-trivial foliation by orbits of a complex Lie subgroup of $G(\mathbb{C})$. This situation is restrictive enough that we are then able to argue that such a foliation can only arise from a Hermitian symmetric foliation on $S$ itself; we give the details in \autoref{hermfolsec}.

\subsection{Acknowledgements}

The author thanks Jacob Tsimerman for numerous helpful discussions and suggestions; in particular, for suggesting that it should be possible to put a $\Qbar$-structure on the normalization of period images, and for introducing the author to jet-theoretic methods in Hodge theory.

\section{Jet Theory}
\label{jetsec}

In this section we develop the theory of higher-dimensional jet spaces, giving an abstract algebraic, an explicit algebraic, and an analytic perspective, all of which will be useful. Much of the work in this section is a straightforward generalization of the usual jet scheme theory, for which the reader may refer to \cite{jetsch3}. 

\subsection{Algebraic Theory}
\label{algjetconstr}

We denote by $\FSch_{K}$ the category of separated schemes of finite type over $K$. Our goal in this section is to construct a family of functors $J^{d}_{r} : \FSch_{K} \to \FSch_{K}$, indexed by non-negative integers $d$ and $r$, and to establish their properties. When we refer to $K$-schemes in this section we will always mean objects of $\FSch_{K}$. We define $A^{d}_{r} = K[t_{1}, \hdots, t_{d}]/(t_{1}, \hdots, t_{d})^{r+1}$, and $\mathbb{D}^{d}_{r} = \Spec A^{d}_{r}$. We view $\mathbb{D}^{d}_{r}$ as a formal $d$-dimensional disk of order $r$. For all $r' \geq r$ and all $d$ we fix a collection of compatible embeddings $\mathbb{D}^{d}_{r} \hookrightarrow \mathbb{D}^{d}_{r'}$ given by $t_{i} \mapsto t_{i}$. From these embeddings we will obtain restriction maps $J^{d}_{r'} \to J^{d}_{r}$.

\begin{defn}
\label{jetspacedef}
Suppose that $S$ is an object in $\FSch_{K}$. We define $J^{d}_{r} S$ to be a scheme representing the contravariant functor $\FSch_{K} \to \Set$ given by 
\[ T \mapsto \Hom_{K}(T \times_{K} \mathbb{D}^{d}_{r}, S), \hspace{1.5em} [T \to T'] \mapsto [\Hom_{K}(T' \times_{K} \mathbb{D}^{d}_{r}, S) \to \Hom_{K}(T \times_{K} \mathbb{D}^{d}_{r}, S)] , \]
where the natural map $\Hom_{K}(T' \times_{K} \mathbb{D}^{d}_{r}, S) \to \Hom_{K}(T \times_{K} \mathbb{D}^{d}_{r}, S)$ obtained by pulling back along $T \times_{K} \mathbb{D}^{d}_{r} \to T' \times_{K} \mathbb{D}^{d}_{r}$. We note that all maps here are of $K$-schemes.
\end{defn}

\begin{prop}
\label{jetreparg}
The functor in \autoref{jetspacedef} is representable.
\end{prop}

\begin{proof}
Call the functor in question $h$, where we leave $r, d$ and $S$ implicit in the notation. By the general theory in \cite{neronmodels}, the Weil restriction $\textrm{Res}_{\mathbb{D}^{d}_{r}/K} (S \times_{K} \mathbb{D}^{d}_{r})$ exists, and is separated and finite type over $K$. It thus suffices to give an isomorphism between $h$ and the functor $h'$ defining the Weil restriction, which is given by 
\[ T \mapsto \Hom_{\mathbb{D}^{d}_{r}}(T \times_{K} \mathbb{D}^{d}_{r}, S \times_{K} \mathbb{D}^{d}_{r}) , \]
on objects, and on maps by pulling back along $T \times_{K} \mathbb{D}^{d}_{r} \to T' \times_{K} \mathbb{D}^{d}_{r}$. Any $\mathbb{D}^{d}_{r}$-map $f : T \times_{K} \mathbb{D}^{d}_{r} \to S \times_{K} \mathbb{D}^{d}_{r}$ is naturally a pair of $K$-maps, one into $S$ and the other into $\mathbb{D}^{d}_{r}$, with the condition that $f$ be a $\mathbb{D}^{d}_{r}$-map translating to the condition that $T \times_{K} \mathbb{D}^{d}_{r} \to \mathbb{D}^{d}_{r}$ be the structure map. Thus we have a bijection $h'(T) \to h(T)$, and it is easy to check that this bijection is natural in $T$.
\end{proof}

\begin{defn}
\label{Jonmaps}
Suppose that $f : S \to S'$ is a map in $\FSch_{K}$. Define $J^{d}_{r} f : J^{d}_{r} S \to J^{d}_{r} S'$ to be the map induced by the natural morphism of functors
\[ [T \mapsto \Hom_{K}(T \times_{K} \mathbb{D}^{d}_{r}, S)] \to [T \mapsto \Hom_{K}(T \times_{K} {\mathbb{D}^{d}_{r}}, S')], \]
induced by post-composition. We observe that this morphism is natural in $T$ as pre- and post-composition commute.
\end{defn}

\autoref{Jonmaps} makes $J^{d}_{r}$ into a functor. Moreover,  for each $r' \geq r$, it is easy to check that the fixed embeddings $\mathbb{D}^{d}_{r} \hookrightarrow \mathbb{D}^{d}_{r'}$ give morphisms of functors $J^{d}_{r'} \to J^{d}_{r}$, and hence in particular give algebraic maps $J^{d}_{r'} S \to J^{d}_{r} S$ which take a jet $j : \mathbb{D}^{d}_{r'} \to S$ to its restriction $\mathbb{D}^{d}_{r} \hookrightarrow \mathbb{D}^{d}_{r'} \xrightarrow{j} S$ to $\mathbb{D}^{d}_{r}$. Moreover, these maps are compatible, in the sense that if $r' \geq r'' \geq r$ then the maps $J^{d}_{r'} \to J^{d}_{r}$ and $J^{d}_{r'} \to J^{d}_{r''} \to J^{d}_{r}$ are equal. As a special case we obtain natural projections $J^{d}_{r} S \to S$ (taking $r = 0$), which associate to a jet $j : \mathbb{D}^{d}_{r} \to S$ the point which it maps to. 

Another construction we will need is that of ``non-degenerate'' jets. This will give, for each $K$-scheme $S$, an open subscheme $J^{d}_{r, nd} S \subset J^{d}_{r} S$ consisting of those jets which induce embeddings of tangent spaces. To give a precise definition of this, we observe that for $r \geq 1$ we have maps $\delta^{d}_{r} : A^{d}_{r} \to K[t_{1}]/( t_{1}^2 ) \times \cdots \times K[t_{d}]/( t_{d}^2 )$ given by $t_{i} \mapsto t_{i}$. As the space representing maps $\Spec K[t]/(t^2) \times T \to S$ is the tangent space $TS$ of $S$, the maps $\delta^{d}_{r}$ induce maps $J^{d}_{r} S \to (TS)^{d}$ given by projecting onto the $d$ tangent vectors $\Spec K[t_{i}]/(t_{i}^2) \to S$ thus induced. We define $J^{d}_{r, nd} S \subset J^{d}_{r} S$ to be the fibre above the open subscheme of $(TS)^{d}$ consisting of those tuples of vectors which span a subspace of $TS$ of dimension $d$.

\subsection{Explicit Models}
\label{exmodelsec}

The abstract theory of jet spaces is useful, however as input to our construction of the maps $\eta^{d}_{r}$ in \autoref{arithdiffcorrespsec} we will need to work with explicit affine models of jet spaces. We therefore describe, in this section, how to explicitly compute jet spaces and the maps between them for subschemes of affine space. This is essentially just a specialization of the Weil restriction construction.

Suppose that $S \subset \mathbb{A}^{n} = \Spec K[z_{1}, \hdots, z_{n}]$ is a closed affine subscheme. Denote by $\mathcal{P}^{d}_{r}$ the set of partitions of integers $0$ through $r$ which have $d$ terms, and let $\ell = |\mathcal{P}^{d}_{r}|$. A map $\mathbb{D}^{d}_{r} \to \mathbb{A}^{n}$ is naturally identified with $n$ formal sums $\sum_{p \in \mathcal{P}^{d}_{r}} a_{p, i} \overline{t}^{p}$, where $1 \leq i \leq n$ and we use multi-index notation to exponentiate the vectors $\overline{t} = (t_{1}, \hdots, t_{d})$. We thus see that $J^{d}_{r} \mathbb{A}^{n}$ is naturally identified with $\mathbb{A}^{n \ell} \simeq \mathbb{A}^{n} \times \mathbb{A}^{n(\ell - 1)}$. Denoting by $\varnothing$ the empty partition, we may also make the identification $a_{\varnothing, i} = z_{i}$.

Now suppose $S \subset \mathbb{A}^{n}$ is cut out by the ideal $I = (f_{1}, \hdots, f_{k})$. Then we may consider, for each $j$, the formal expansion
\begin{equation} 
\label{affjetspaceeq}
f_{j}\left( \sum_{p \in \mathcal{P}} a_{p, 1} \overline{t}^{p}, \hdots, \sum_{p \in \mathcal{P}} a_{p, n} \overline{t}^{p} \right) , 
\end{equation}
which is a polynomial in the $\overline{t}^{p}$. We wish to impose the condition that this formal expression is zero up to order $r$, which gives a total of $\ell$ polynomial equations in the coordinates $a_{p, i}$ determined by the $\ell$ coefficients of the terms $\overline{t}^{p}$ with order less than or equal to $r$. Taking all $f_{j}$ into consideration for $1 \leq j \leq k$, we obtain a total of $\ell k$ polynomial equations in the coordinates on $J^{d}_{r} \mathbb{A}^{n}$. We define $J^{d}_{r} S$ to be the $K$-scheme cut out by these equations.

Now suppose we have two closed affine subschemes $S \subset \mathbb{A}^{n}$ and $S' \subset \mathbb{A}^{m}$, an algebraic map $g : S \to S'$, and we construct $J^{d}_{r} S \subset J^{d}_{r} \mathbb{A}^{n}$ and $J^{d}_{r} S' \subset J^{d}_{r} \mathbb{A}^{m}$ as above. Then in these coordinates, $g$ has $m$ component functions $g_{j}$ for $1 \leq j \leq m$. If $j : \mathbb{D}^{d}_{r} \to S \subset \mathbb{A}^{n}$ is an element of $J^{d}_{r} S$, we may compute the composition $g \circ j$ by evaluating the formal expressions
\begin{equation}
\label{affjetmapeq}
g_{j}\left(\sum_{p \in \mathcal{P}} a_{p, 1} \overline{t}^{p}, \hdots, \sum_{p \in \mathcal{P}} a_{p, n} \overline{t}^{p}\right) , 
\end{equation}
and truncating after order $r$. The coefficients of these formal expressions once again define polynomials in the $a_{p, i}$, which give us our maps $J^{d}_{r} S \to J^{d}_{r} S'$.

Both the defining equations for $J^{d}_{r} S$ and polynomials giving the map $J^{d}_{r} S \to J^{d}_{r} S'$ can also be computed in a different way, which will be useful in both the analytic theory that follows and our constructions in \autoref{jetdiffeq}. To compute the coefficient $c_{p}$ of the $\overline{t}^{p}$ term in \autoref{affjetspaceeq}, for instance, we may use several applications of the multivariate chain rule so that the constant term of the resulting formal expression is a positive integer multiple of $c_{p}$. Dividing out this integer, we find that there exist $\mathbb{Q}$-polynomials in the $a_{p, i}$ and the partial derivatives $\partial_{q} f_{j} (a_{\varnothing, 1}, \hdots, a_{\varnothing, n})$ which define $J^{d}_{r} S \subset J^{d}_{r} \mathbb{A}^{n}$, where $\partial_{q}$ denotes a sequence of up to $r$ partial derivatives with respect to the coordinates $a_{\varnothing, i}$. These polynomials are also universal, in the sense recorded by the following proposition:

\begin{notn}
Define $\mathcal{Q}_{r} = \bigcup_{i = 0}^{r} \{ 1, \hdots, n \}^{i}$. Then for $q = (q_{1}, \hdots, q_{k_{i}}) \in \mathcal{Q}_{r}$ we let $\partial_{q} = \partial_{q_{k_{i}}} \cdots \partial_{q_{1}}$.
\end{notn}

\begin{prop}
\label{univjeteq}
Let $R$ be the $\mathbb{Q}$-algebra freely generated by the formal symbols $a_{p, i}$ for $1 \leq i \leq n$ and $p \in \mathcal{P}^{d}_{r}$, and $\partial_{q} f_{j} (a_{\varnothing, 1}, \hdots, a_{\varnothing, n})$ for $1 \leq j \leq k$ and $q \in \mathcal{Q}_{r}$; let $\ell = |\mathcal{P}^{d}_{r}|$. Then there exists a map $r : \Spec R \to \mathbb{A}^{k \ell}$ that satisfies the following properties:

\begin{itemize}
\item[(i)] The map $r$ is $\mathbb{Q}$-algebraic.
\item[(ii)] Let $f = (f_{1}, \hdots, f_{k})$ be functions on $\mathbb{A}^{n}$ defining a subscheme $S$, and let $\pi : J^{d}_{r} \mathbb{A}^{n} \to \mathbb{A}^{n}$ be the natural projection. Consider the map $\varphi_{f} : J^{d}_{r} \mathbb{A}^{n} \to \Spec R$ defined by 
\[ \sigma \mapsto (\sigma, \partial_{q} f_{j} (\pi(\sigma)))  \]
where $q$ ranges over all elements of $\mathcal{Q}_{r}$ and $j$ ranges from $1$ to $k$. Then the subscheme $J^{d}_{r} S \subset J^{d}_{r} \mathbb{A}^{n}$ is the fibre above $0$ of $r \circ \varphi_{f}$.
\end{itemize}
\end{prop}

\begin{proof}
This is just a rephrasing of what we have said. Indeed, the subscheme $J^{d}_{r} S \subset J^{d}_{r} \mathbb{A}^{n}$ is defined by $\ell k$ equations obtained from expanding \autoref{affjetspaceeq} for $1 \leq j \leq k$ and setting the coefficient of $\overline{t}^{p}$ to zero for all $p$. One may calculate these coefficients for indeterminate choices of $f_{j}$ by taking partial derivatives of the expression in \autoref{affjetspaceeq}, applying the multivariate chain rule, and possibly dividing out by an integer coming from the exponents of $\overline{t}^{p}$, and evaluating at $\overline{t} = 0$. This gives a collection of $\ell k$ functions which are $\mathbb{Q}$-algebraic in the coordinates of $\Spec R$, and defines the map $r$. That $J^{d}_{r} S \subset J^{d}_{r} \mathbb{A}^{n}$ is the fibre above $r \circ \varphi_{f}$ is then nothing but the statement that the common vanishing locus of these functions, after replacing the indeterminate derivatives of the $f_{j}$ with the functions themselves, define the locus $J^{d}_{r} S$.
\end{proof}

The analogous fact is true for maps, which we also record.

\begin{prop}
\label{univjetmapeq}
Let $\Spec R$ be as in \autoref{univjeteq}, except with $g_{j}$ replacing $f_{j}$ and $m$ replacing $k$; let $\ell = |\mathcal{P}^{d}_{r}|$. Then there exists a map $r : \Spec R \to J^{d}_{r} \mathbb{A}^{m}$ which satisfies the following properties:
\begin{itemize}
\item[(i)] The map $r$ is $\mathbb{Q}$-algebraic.
\item[(ii)] Let $S \subset \mathbb{A}^{n}$ and $S' \subset \mathbb{A}^{m}$ be any two closed subschemes, and let $g : S \to S'$ be a map between them. View $g$ as a map $\mathbb{A}^{n} \to \mathbb{A}^{m}$ with components $g_{j}$ for $1 \leq j \leq m$, and define $\varphi_{g} : J^{d}_{r} \mathbb{A}^{n} \to \Spec R$ by
\[ \sigma \mapsto (\sigma, \partial_{q} g_{j}(\pi(\sigma))) \]
where $q$ ranges over $\mathcal{Q}_{r}$ and $j$ ranges from $1$ to $m$. Then the induced map $J^{d}_{r} g : J^{d}_{r} S \to J^{d}_{r} S'$ is the restriction of $r \circ \varphi_{g}$ to $J^{d}_{r} S \subset J^{d}_{r} \mathbb{A}^{n}$.
\end{itemize}
\end{prop}

\begin{proof}
The map $r$ is defined in the same way as in \autoref{univjeteq} except using \autoref{affjetmapeq} instead of \autoref{affjetspaceeq}. The second property is also immediate for the same reason. 
\end{proof}

\subsection{Analytic Theory}
\label{anjetspace}

Denote by $\An$ the category of complex analytic spaces. In this section we describe the analytic analogue of the construction in the previous section, which defines a functor $J^{d}_{r} : \An \to \An$ which associates to an analytic space $X$ an analytic space $J^{d}_{r} X$. We note that $\mathbb{D}^{d}_{r}$ is also naturally an analytic space, where we take $K = \mathbb{C}$. 

The majority of the work is done by \autoref{exmodelsec}. Indeed, it is clear that it suffices to give the construction of $J^{d}_{r} X$ when $X$ is a basic analytic space. Thus, we may regard $X$ as a closed analytic subvariety of an open set $U \subset \mathbb{C}^{n}$ defined by analytic functions $f_{1}, \hdots, f_{k}$ on $U$. We define $J^{d}_{r} X$ as an analytic subspace of an open set in $\mathbb{C}^{n\ell}$ where $\ell = |\mathcal{P}^{d}_{r}|$ and the coordinates on $\mathbb{C}^{n \ell}$ are given by $a_{i, p}$ as before. To avoid manipulations with power series, we may simply define $J^{d}_{r} X$ to be given by the fibre above $0$ of $\an{r} \circ \varphi_{f}$, where $\varphi_{f}$ and $r$ are as in \autoref{univjeteq}. It is then clear that if $S$ is a separated finite-type $\mathbb{C}$-scheme, then $\an{(J^{d}_{r} S)} = J^{d}_{r} \an{S}$, as this fact reduces to the affine case and the affine case is automatic.

Given a map $g : X \to X'$ of analytic spaces, we also wish to describe the induced map $J^{d}_{r} g : J^{d}_{r} X \to J^{d}_{r} X'$. Once again it suffices to handle the case where $X$ and $X'$ are basic analytic spaces. We therefore may assume that $X$ is a closed analytic subvariety of an open set $U \subset \mathbb{C}^{n}$, and that $X'$ is a closed analytic subvariety of an open set $U' \subset \mathbb{C}^{m}$. In these coordinates, the map $g$ has components $g_{j}$ for $1 \leq j \leq n$. We define the maps $J^{d}_{r} g$ as the composition $\an{r} \circ \varphi_{g}$, where $r$ and $\varphi_{g}$ are as in \autoref{univjetmapeq}, except $\varphi_{g}$ is only defined on $J^{d}_{r} U$ instead of all of $J^{d}_{r} \mathbb{C}^{n}$. The functor $J^{d}_{r}$ is then extended to all maps by gluing, and it is once again clear that the construction is compatible with analytification.

We note that associated to the embeddings $\mathbb{D}^{d}_{r} \hookrightarrow \mathbb{D}^{d}_{r'}$ we have restriction functors as in the preceding section; we omit the verification. By the same construction as in the \autoref{algjetconstr}, we have open analytic subspaces $J^{d}_{r, nd} X \subset J^{d}_{r} X$ consisting of those maps $\mathbb{D}^{d}_{r} \hookrightarrow X$ which are injective on tangent spaces. We now give some applications of non-degenerate jets that we will find useful.

\begin{lem}
\label{analndissur}
Suppose that $X$ is an analytic space, and that $j \in J^{d}_{r, nd} X$. Then the map $j : \mathbb{D}^{d}_{r} \to X$ is an embedding, in the sense that if $x \in X$ is the basepoint of $j$, then the induced map $j^{\sharp} : \mathcal{O}_{X, x} \to A^{d}_{r}$ is surjective.
\end{lem} 

\begin{proof}
To show that the induced map $j^{\sharp}: \mathcal{O}_{S, s} \to A^{d}_{r}$ is surjective, it suffices to show that the generators $t_{1}, \hdots, t_{d}$ are in the image. We know that after composing with the map $\delta^{d}_{r}$ from \autoref{algjetconstr} we obtain a linearly independent set of vectors, and hence the map $\delta^{d}_{r} \circ j^{\sharp}$ must be surjective. It follows that there are elements $z_{1}, \hdots, z_{d} \in \mathcal{O}_{S, s}$ such that $(\delta^{d}_{r} \circ j^{\sharp})(z_{i}) = \delta^{d}_{r}(j^{\sharp}(z_{i})) = t_{i}$, and so in particular $j^{\sharp}(z_{i})$ lies in the coset $t_{i} + (t_{i}^2, \{ t_{j} : j \neq i \})$. It then follows by explicit calculation that $(j^{\sharp}(z_{1}), \hdots, j^{\sharp}(z_{d})) = (t_{1}, \hdots, t_{d})$, and so surjectivity follows.
\end{proof}

\begin{prop}
\label{jetdimlem}
Suppose that $X$ is an analytic space, that $x \in X$ is a point, and the fibre of $J^{d}_{r, nd} X$ above $x \in X$ is non-empty for all $r$. Then $X$ has dimension at least $d$ at $x$.
\end{prop}

\begin{proof}
Let $m$ be the dimension of the noetherian local ring $\mathcal{O}_{X, x}$. We will show that $m \geq d$. By  \cite[\href{https://stacks.math.columbia.edu/tag/00KQ}{Tag 00KQ}]{stacks-project}, the statement that $\mathcal{O}_{X, x}$ has dimension $m$ is equivalent to the statement that there exists an ideal $I \subset \mathcal{O}_{X, x}$ such that $\sqrt{I} = \mathfrak{m}_{X, x}$, and $I$ is generated by $m$ elements $z_{1}, \hdots, z_{m}$. As $\mathcal{O}_{X, x}$ is a Noetherian ring, by \cite[\href{https://stacks.math.columbia.edu/tag/00IM}{Tag 00IM}]{stacks-project} we eventually have that $\mathfrak{m}_{X, x}^{k} \subset I$ for some $k$. 

Suppose $j$ is a non-degenerate jet in the fibre of $J^{d}_{r, nd} X$ above $x$, suppose $r \geq k$, and consider the induced morphism $j^{\sharp} : \mathcal{O}_{X, x} \to A^{d}_{r}$. This map is surjective by \autoref{analndissur}. Since $r \geq k$, we must have $\textrm{ker}(j^{\sharp}) \subset \mathfrak{m}_{X, x}^{k}$, since surjectivity implies that the generators of $\mathfrak{m}_{X, x}$ must map to generators of the maximal ideal in $A^{d}_{r}$. As we have both that $j^{\sharp}$ is surjective and $\textrm{ker}(j^{\sharp}) \subset I$, an elementary exercise in commutative algebra shows that $j^{\sharp}(\sqrt{I}) = \sqrt{j^{\sharp}(I)}$. We thus have that $\sqrt{j^{\sharp}(I)} = ( t_{1}, \hdots, t_{d} )$, which implies that the minimal number of generators of $j^{\sharp}(I)$, and hence $I$, is at least $d$. Thus $m \geq d$. 
\end{proof}

\begin{defn}
\label{compjetdef}
Given two jets $j_{r} \in J^{d}_{r} X$ and $j_{r'} \in J^{d}_{r'} X$, we say that $j_{r'}$ and $j_{r}$ are \emph{compatible} if $\pi^{r'}_{r}(j_{r'}) = j_{r}$. Given an infinite sequence $\{ j_{r} \}_{r \geq 0}$ with $j_{r} \in J^{d}_{r} X$, we will often refer to the $j_{r}$ as a \emph{compatible family} or \emph{compatible sequence} if $j_{r}$ is compatible with $j_{r+1}$ for all $r$.
\end{defn}

\begin{defn}
\label{imdef}
Given a map $f : (X, x) \to (Y, y)$ of analytic germs, we call the smallest $(Z, y) \subset (Y, y)$ through which $f$ factors the image of $f$.
\end{defn}

\begin{prop}
\label{sameimagelem}
Let $X_{1}, X_{2}$ and $Y$ be smooth analytic spaces, and let $\alpha_{i} : X_{i} \to Y$ be two analytic maps. Suppose that $d = \dim X_{i}$ for $i = 1, 2$, and $j_{r, i} \in J^{d}_{r,nd} X$ for $r \geq 0$ and $i = 1, 2$ is an infinite family of compatible non-degenerate jets such that $\alpha_{1} \circ j_{r, 1} = \alpha_{2} \circ j_{r, 2}$ for all $r$. Then if we let $x_{i} = j_{0,i}$ and $y = \alpha_{i}(x_{i})$ the maps $\alpha_{i} : (X_{i}, x_{i}) \to (Y, y)$ of analytic germs have the same image.
\end{prop}

\begin{proof}
The statement is local, so we assume that $X_{i} = \mathbb{C}^{d}$ and $Y = \mathbb{C}^{m}$. Denote by $A^{d}_{\infty}$ the projective limit of the $A^{d}_{r}$; this is simply the ring of formal power series in $d$ variates. We then denote by $\mathcal{O}_{d}$ and $\mathcal{O}_{m}$ the rings of germs of holomorphic functions at the origin in $\mathbb{C}^{d}$ and $\mathbb{C}^{m}$ respectively. Then the family of compatible jets $j_{r, i}$ for $i = 1, 2$ induces a single morphism $j_{\infty, i} : \mathcal{O}_{d} \to A^{d}_{\infty}$, and we are reduced to the following statement: given two morphisms $\eta_{i} : \mathcal{O}_{m} \to \mathcal{O}_{d}$ for $i = 1, 2$ with the property that $j_{\infty, 1} \circ \eta_{1} = j_{\infty, 2} \circ \eta_{2}$, then the maps $\eta_{i}$ both have the same kernel. To prove this statement we are free to replace the maps $\eta_{i}$ with their formal completions $\widehat{\eta}_{i}$ which are maps $\widehat{\mathcal{O}}_{m} \to \widehat{\mathcal{O}}_{d}$. Then as the jets $j_{r}$ are non-degenerate, the maps $\widehat{j}_{\infty, i} : \widehat{\mathcal{O}}_{d} \to A^{d}_{\infty}$ are isomorphisms by the formal inverse function theorem, so the result follows.
\end{proof}

\section{Arithmetic Differential Correspondences}
\label{arithdiffcorrespsec}

In this section we construct $K$-algebraic maps $J^{d}_{r} S \to G \backslash J^{d}_{r} \ch{D}$ associated to $K$-algebraically-defined polarized variations of Hodge structure. We will often implicitly fix a choice of $d$ and $r$, as for instance in the previous sentence, with the understanding that the claims we make apply for all choices of $d$ and $r$.

\subsection{Jets Defined by Linear Differential Equations}
\label{jetdiffeq}

Suppose that $U$ is a smooth $K$-variety, and fix an integer $m$. In this section we are interested in constructing certain $K$-algebraic \emph{jet evaluation maps} associated to a system of differential equations of the form 
\begin{equation}
\label{defdiffeq3}
 df_{jk} = \sum_{i = 1} f_{ik} c_{ij} ,
\end{equation}
where $1 \leq i, j, k \leq m$ and the $c_{ij}$ are global sections of $\Omega^{1}_{U}$. We will only be interested in solutions such that the matrix $[f_{ij}]$ is invertible; we think of these solutions as being analytic maps $f : B \to \GL_{m}$, where $\GL_{m}$ denotes the space of invertible $m \times m$ matrices, and $B \subset U$ is an analytic neighbourhood. We understand the equality of \autoref{defdiffeq3} as taking place inside $\Omega^{1}_{\an{U}}(B)$. 

We will see in \autoref{alphamapconstr} that a system of differential equations like that of \autoref{defdiffeq3} arises naturally when computing a basis of flat sections of a variation of Hodge structures. The fact that solutions always exist and that their germs at a point $s \in U$ are uniquely determined by an initial condition $f(s) = M$ will be automatic in the variation of Hodge structures case, so we assume it in the following proposition.

\begin{prop}
\label{keymapprop}
Suppose that for any $s \in U$, solutions to \autoref{defdiffeq3} exist on some analytic neighbourhood $B \subset U$ of $s$ and are uniquely determined by the initial condition $f(s) = M$. Then there exists a $K$-algebraic map $\beta : J^{d}_{r} U \times \GL_{m} \to J^{d}_{r} \GL_{m}$, which has the property that if $f : B \to \GL_{m}$ is an analytic solution to the differential system of \autoref{defdiffeq3}, then $\beta(\sigma, f(\pi(\sigma))) = f \circ \sigma$ for all $\sigma \in J^{d}_{r} U$, where $\pi : J^{d}_{r} U \to U$ is the usual projection.
\end{prop}

\begin{proof}
We first claim that it suffices to do this construction Zariski-locally on $U$. Indeed, suppose that $U_{1}$ and $U_{2}$ are two $K$-algebraic affine neighbourhoods of $U$, and we have maps $\beta_{i} : J^{d}_{r} U_{i} \times \GL_{m} \to J^{d}_{r} \GL_{m}$ satisfying the stated properties with respect to the restrictions of \autoref{defdiffeq3} to $U_{1}$ and $U_{2}$. As we are dealing with maps of varieties, it suffices to check that the maps agree on complex points. Letting $(\sigma, M) \in J^{d}_{r} (U_{1} \cap U_{2}) \times \GL_{m}$ be a complex point, we may construct a solution $f : B \to \GL_{m}$ to the differential system on a neighbourhood $B \subset U_{1} \cap U_{2}$ such that $f(s) = M$, where $s$ is the basepoint of $\sigma$. Then we see that
\[ \beta_{1}(\sigma, M) = f \circ \sigma = \beta_{2}(\sigma, M) , \]
so we see that the two maps agree.

As $U$ is smooth, we can now restrict to the case where the sheaf $\Omega^{1}_{U}$ is free and trivialized by sections $dz_{1}, \hdots, dz_{n}$. Such a trivialization allows us to define \emph{partial derivatives} with respect to the $z_{i}$ for any function $g$ on $U$: we write $dg = \sum_{i = 1}^{n} a_{i} d z_{i}$ for unique functions $a_{i}$ and define $\partial_{\ell} g = a_{\ell}$. The differential equation \autoref{defdiffeq3} then becomes
\begin{equation}
\label{defdiffeq2}
 \partial_{\ell} f_{jk} = \sum_{i} f_{ik} c_{ij, \ell} 
\end{equation}
for algebraic functions $c_{ij, \ell}$ on $U$. This trivialization of $\Omega^{1}_{U}$ has associated to it an \'etale map $\rho : U \to \mathbb{A}^{n} = \Spec K[x_{1}, \hdots, x_{n}]$ defined by $x_{i} \mapsto z_{i}$. We now make two claims which are crucial to our construction:
\begin{claim}
\label{claim1}
Let $g : B \to \mathbb{C}$ be an analytic function defined on some open subset $B \subset U$ such that $\rho$ restricts to an isomorphism on $B$. Define $g' = g \circ \rho^{-1}$, and consider partial derivatives of $g'$ with respect to the coordinates $x_{1}, \hdots, x_{n}$. Then if $q \in \mathcal{Q}_{r}$ (see the definition preceding \autoref{univjeteq} and \autoref{univjetmapeq}), then we have
\[ (\partial_{q} g')(x) = (\partial_{q} g)(\rho^{-1}(x)) \]
for any $x \in \rho(B)$, where the partial derivatives of $g$ are taken with respect to $z_{1}, \hdots, z_{n}$.
\end{claim}

\begin{proof}
We note that the claim holds by definition in the case where $q = \varnothing$, so we may argue by induction; we thus need simply to show the statement for the case when $\partial_{q} = \partial_{\ell}$ is a single partial derivative. We have that
\[ d(g \circ \rho^{-1})_{x} = (dg)_{\rho^{-1}(x)} \circ (dp)^{-1}_{x} . \]
The map $(d\rho)^{-1}$ sends the sections $\frac{\partial}{\partial x_{i}}$ to $\frac{\partial}{\partial z_{i}}$, so we find that
\[ (\partial_{\ell} g')(x) = d(g \circ p^{-1})_{x}\left( \frac{\partial}{\partial x_{i}} \right) = (d g)_{p^{-1}(x)} \left( \frac{\partial}{\partial z_{i}} \right) = (\partial_{\ell} g)(p^{-1}(x)) . \]
\end{proof}

\begin{claim}
\label{claim2}
Suppose that $f : B \to \GL_{m}$ is a solution to the differential system of \autoref{defdiffeq2}. Then for any $q \in \mathcal{Q}_{r}$ and any component $f_{jk}$ of $f$ there is a polynomial $\xi_{q, jk}$ in the components of $f$ whose coefficients are algebraic functions on $U$ (so if $U = \Spec A$, then $\xi_{p, jk}$ is an element of $A[f_{tu}, 1 \leq t, u \leq m]$) such that
\[ \partial_{q} f_{jk} = \xi_{q, jk}([f_{tu}]) \]
is an equality of functions on $B$. Moreover, this polynomial is independent of $f$.
\end{claim}

\begin{proof}
This claim is trivial when $q = \varnothing$, and follows when $\partial_{q} = \partial_{\ell}$ is a single partial derivative by \autoref{defdiffeq2}. The case of higher-order partial derivatives is obtained by repeatedly differentiating \autoref{defdiffeq2} and using the polynomials for the lower-order cases to obtain the higher-order cases by substitution.
\end{proof}

We now give an informal description of how $\beta$ is constructed. For any function $f : B \to \GL_{m}$ and any element $\sigma \in J^{d}_{r} U$, the composition $f \circ \sigma$ is, when written in coordinates, a $\mathbb{Q}$-algebraic function of the coordinates of $\sigma$ and the partial derivatives of $f$. The content of \autoref{claim2} can be interpreted as saying that the partial derivatives $\partial_{q} f_{jk}(\pi(\sigma))$ are $K$-algebraic functions of $\pi(\sigma)$ and $f(\pi(\sigma))$. Thus, by expressing  the derivatives $\partial_{q} f_{jk}(\pi(\sigma))$ algebraically over $K$ in terms of $f(\pi(\sigma))$ and $\pi(\sigma)$, we find that $f \circ \sigma$ is a $K$-algebraic function of $\sigma$ and $f(\pi(\sigma))$.

Proceeding now more formally, let $\mathbb{M}$ be the space of $m \times m$ matrices, of which $\GL_{m}$ is a natural subspace, and let $r : \Spec R \to J^{d}_{r} \mathbb{M}$ be the map of \autoref{univjetmapeq} associated to the affine spaces $\mathbb{A}^{n}$ and $\mathbb{M}$. The map $r$ is a $\mathbb{Q}$-algebraic map with the following property: if $g : B \to \mathbb{M}$ is any analytic map on $B \subset \mathbb{A}^{n}$ and with components $g_{ij}$, and we define $\varphi_{g} : J^{d}_{r} B \to \an{(\Spec R)}$ via
\[ \sigma \mapsto (\sigma, \partial_{q} g_{ij}(\pi(\sigma))) , \]
where $q$ ranges over $\mathcal{Q}_{r}$ and $ij$ ranges over all possible indices, then $g \circ \sigma = (\an{r} \circ \varphi_{g})(\sigma)$ for all $\sigma \in J^{d}_{r} B$. 

Denote by $a_{tu}$ the coordinates on $\mathbb{M}$. We note that if $U = \Spec A$, then the coordinate ring of $U \times \mathbb{M}$ is $A[a_{tu}, 1 \leq t, u \leq m]$, so in particular the functions $\xi_{q, jk}$ of \autoref{claim2} may be interpreted as algebraic functions on $U \times \GL_{m}$. We now define $\zeta : J^{d}_{r} U \times \GL_{m} \to \Spec R$ to be the $K$-algebraic map
\[ (\sigma, [a_{tu}]) \mapsto (\rho \circ \sigma, \xi_{q,jk}([a_{tu}], \pi(\sigma))) , \]
where $q$ ranges over all elements of $\mathcal{Q}_{r}$ and $1 \leq j, k \leq m$. We define $\beta = r \circ \zeta$. It now suffices to consider an analytic function $f : B \to \GL_{m}$ defined on $U$ which satisfies \autoref{defdiffeq2}, and check that $\beta(\sigma, f(\pi(\sigma))) = f \circ \sigma$. We may assume that $\rho$ restricts to an isomorphism on $B$, and let $f' = f \circ \rho^{-1}$. We have that
\begin{align*}
f \circ \sigma &= f' \circ \rho \circ \sigma \\
&= (\an{r} \circ \varphi_{f'})(\rho \circ \sigma) \\
&= \an{r}(\rho \circ \sigma, (\partial_{q} f'_{jk})(\pi(\rho \circ \sigma))) \\
&= \an{r}(\rho \circ \sigma, (\partial_{q} f_{jk})(\pi(\sigma))) & (\star) \\
&= \an{r}(\rho \circ \sigma, \xi_{q, ij}(f_{jk}(\pi(\sigma)), \pi(\sigma))) \\
&= \an{r}(\zeta(\sigma, [f_{jk}(\pi(\sigma))])) \\
&= \beta(\sigma, f(\pi(\sigma))) ,
\end{align*}
where on the line labelled $(\star)$ we have applied \autoref{claim1}.
\end{proof}

\begin{lem}
\label{rightact}
The map of \autoref{keymapprop} is invariant under the natural right-action of $\GL_{m}$.
\end{lem}

\begin{proof}
To check the invariance of $\beta$ under the right-action it suffices to check that if $f$ solves the differential system of \autoref{defdiffeq3}, then so does $f A$, where $A$ is any matrix in $\GL_{m}$. Then
\[ \beta(\sigma, M A) = (f A) \circ \sigma = (f \circ \sigma) A = \beta(\sigma, M) A . \]
\end{proof}

\subsection{Jets Defined by Local Period Maps}
\label{alphamapconstr}

\subsubsection{Differential Equations for Flat Frames}

We now apply the construction in the preceding section to the setting of a variation of Hodge structure, where a differential system of the form required in \autoref{keymapprop} is supplied by the connection $\nabla$. To demonstrate how this works, we suppose that $\mathbb{V}$ is a variation of Hodge structure on a smooth $K$-variety $S$ with Hodge bundle $\mathcal{H}$, and that $\nabla : \mathcal{H} \to \Omega^{1}_{S} \otimes \mathcal{H}$ is the Gauss-Manin connection, also defined over $K$. Let $U \subset S$ be a Zariski open neighbourhood on which $\mathcal{H}$ and the filtered pieces $F^{k}$ are trivial. We define
\begin{defn}
A \emph{filtration-compatible frame} $v^{i} : U \to \mathcal{H}$ for $1 \leq i \leq m$ is a $K$-algebraic basis of global sections of $\restr{\mathcal{H}}{U}$ such that some initial segment of the sequence $v^{1}, \hdots, v^{m}$ spans $F^{k}$ for each $k$.
\end{defn}

\noindent Given a filtration-compatible frame $v^{i} : U \to \mathcal{H}$ we write $\nabla v^{i} = \sum_{j = 1}^{m} c_{ij} \otimes v^{j}$ for $c_{ij} \in \Omega^{1}_{U}$. On sufficiently small analytic neighbourhoods $B \subset U$, we may ask for flat sections $b^{k} = \sum_{i = 1}^{m} f_{ik} v^{i}$ in terms of the $v^{i}$, where the $f_{ij}$ are analytic functions on $B$. We then compute that
\begin{align*}
\nabla b^{k} &= \nabla \left(\sum_{i = 1}^{m} f_{ik} v^{i} \right) \\
&= \sum_{j = 1}^{m} d f_{jk} \otimes v^{j} + \sum_{i = 1}^{m} f_{ik} \left( \sum_{j = 1}^{m} c_{ij} \otimes v^{j} \right) \\
&= \sum_{j = 1}^{m} \left( df_{jk} + \sum_{i = 1}^{m} f_{ik} c_{ij} \right) \otimes v^{j} ,
\end{align*}
and so the condition $\nabla b^{k} = 0$ induces a system of differential equations $df_{jk} = - \sum_{i = 1}^{m} f_{ik} c_{ij}$ which agree with \autoref{defdiffeq3} after absorbing the sign into the $c_{ij}$. We note that the existence of such solutions is guaranteed by the existence of the local system $\mathbb{V}_{\mathbb{C}}$ on $S$, and the fact the solutions are determined by the initial condition $f(s_{0}) = M$ at any $s_{0} \in B$ amounts to the fact that $\mathbb{V}_{\mathbb{C}}(B) \to \mathbb{V}_{\mathbb{C}, s_{0}}$ is an isomorphism. (That is, it amounts to the fact that choosing the basis $b^{1}(s_{0}), \hdots, b^{m}(s_{0})$ determines the basis in some neighbourhood of $s_{0}$.) 

We thus obtain a map $\beta: J^{d}_{r} U \times \GL_{m} \to J^{d}_{r} \GL_{m}$ from \autoref{keymapprop} on any open $U \subset S$ on which $\mathcal{H}$ and the filtered pieces $F^{k}$ are trivial. The result of \autoref{keymapprop} now specializes to say 

\begin{prop}
\label{keymapprop2}
Let $\mathbb{V}$ be a $K$-algebraically-defined variation of Hodge structure on $S$, and let $v^{i} : U \to \mathcal{H}$ be a filtration-compatible frame on some $K$-algebraic open subset $U \subset S$. Then there exists a $K$-algebraic map $\beta : J^{d}_{r} U \times \GL_{m} \to J^{d}_{r} \GL_{m}$, which has the property that if $f : B \to \GL_{m}$ is an analytic function such that $b^{k} = \sum_{i = 1}^{m} f_{ik} v^{i}$ gives a flat basis, then $\beta(\sigma, f(\pi(\sigma))) = f \circ \sigma$ for all $\sigma \in J^{d}_{r} U$, where $\pi : J^{d}_{r} U \to U$ is the usual projection. \qed
\end{prop}

\subsubsection{Local Period Maps Revisited}

Let us now revisit the notion of a local period map introduced in the introduction, and use them to construct a refinement of the maps $\beta$ from \autoref{keymapprop2}. We suppose we have a $K$-algebraically-defined polarized variation of Hodge structure $\mathbb{V}$ on the smooth $K$-variety $S$, with fibres of dimension $m$. We fix a polarization form $Q : \mathbb{Z}^{m} \otimes \mathbb{Z}^{m} \to \mathbb{Z}$ of the same type as the polarization $\mathcal{Q}$ on $\mathbb{V}$. We denote by $\ch{L}$ the flag variety parametrizing flags on $\mathbb{C}^{m}$ with the same Hodge numbers as the filtration on $\mathbb{V}$, and we denote by $\ch{D} \subset \ch{L}$ the closed subvariety of those Hodge flags satisfying the first Hodge-Riemann bilinear relation with respect to $Q$. We let $G = \textrm{Aut}(\mathbb{Z}^{m}, Q)$ be the subgroup of $\GL_{m}$ consisting of those elements preserving $Q$. Denoting by $D \subset \ch{D}$ the open subspace parametrizing Hodge flags satisfying the additional positivity condition, the variation $\mathbb{V}$ induces a canonical period map $\varphi : \an{S} \to \Gamma \backslash D$, where $\Gamma = G(\mathbb{Z})$.

There is a natural $\mathbb{Q}$-algebraic quotient map $q : \GL_{m} \to \ch{L}$, invariant under the left action of $\GL_{m}$, which sends a matrix $M \in GL_{m}(\mathbb{C})$ to the Hodge flag $F^{\bullet}$ whose $F^{k}$-subspace is spanned by the first $\dim F^{k}$ columns of $M$. We then make the following claim:

\medskip

\begin{notn}
We denote by $\iota : \GL_{m} \to \GL_{m}$ the inversion map, and by $e^{1}, \hdots, e^{m}$ the standard basis of $\mathbb{Z}^{m}$.
\end{notn}

\begin{lem}
\label{locpermapredef}
With the above setup, a local period map $\psi : B \to \ch{D}$ is exactly a map of the form $q \circ \iota \circ f$, where $f : B \to \GL_{m}$ defines a flat basis $b^{k} = \sum_{i = 1}^{m} f_{ik} v^{i}$ with respect to a filtration-compatible frame $v^{i} : U \to \mathcal{H}$, and $\mathcal{Q}(b^{i}, b^{j}) = Q(e^{i}, e^{j})$ for all $1 \leq i, j \leq m$..
\end{lem}

\begin{proof}
Suppose the basis $b^{k}$ is chosen to be an integral basis of $\mathbb{V}(B)$, whose polarization matrix $\mathcal{Q}(b^{i}, b^{j})$ agrees with $Q(e^{i}, e^{j})$, where $e^{1}, \hdots, e^{m}$ is the standard basis of $\mathbb{Z}^{m}$. Then the entries of the matrix $f^{-1}(s)$ for $s \in B$ are the coordinates of the filtration-compatible frame $v^{1}(s), \hdots, v^{m}(s)$ representing a polarized Hodge flag with respect to an integral basis, hence $q \circ \iota \circ f$ gives a local lift of the canonical period map $\varphi : \an{S} \to \Gamma \backslash D$. 

We have defined local period maps to be maps of the form $g \cdot \psi$ where $\psi$ is a local lift and $g \in G(\mathbb{C})$, so it suffices to show that replacing $b^{k}$ with a more general complex basis $b'^{\ell}$ of $\mathbb{V}_{\mathbb{C}}(B)$ (still with $\mathcal{Q}(b'^{i}, b'^{j}) = Q(e^{i}, e^{j})$) results in a map of the form $g \cdot \psi$ for $g \in G(\mathbb{C})$, and that every map of the form $g \cdot \psi$ is obtained in this way. As we saw in \autoref{rightact} that $\GL_{m}$ acts (necessarily freely and transitively) on functions $f' : B \to \GL_{m}$ defining flat frames in terms of the $v^{i}$, it suffices to show that the $G(\mathbb{C})$-orbit containing $f$ consists of exactly those flat frames $b'^{\ell}$ satisfying $\mathcal{Q}(b'^{i}, b'^{j}) = Q(e^{i}, e^{j})$ for $1 \leq i, j \leq m$. But given a complex basis $b'^{\ell}$ related to $b^{k}$ by $b'^{\ell} = \sum_{r = 1}^{m} g_{r \ell} b^{r}$ for some matrix $g = [g_{r \ell}]$, we have 
\[ g \in G(\mathbb{C}) \hspace{2em} \iff \hspace{2em} \mathcal{Q}(b'^{i}, b'^{j}) = Q(e^{i}, e^{j}) \textrm{ for all }1 \leq i, j \leq m ,\]
hence the result.
\end{proof}

We now construct a refinement of the maps $\beta$ of \autoref{keymapprop2} which can be used to compute the action of $J^{d}_{r} \psi$ whenever $\psi$ is a local period map associated to a filtration-compatible frame $v^{i} : U \to \mathcal{H}$. To ensure we obtain just those maps $\psi$ which respect the polarization, we define the $K$-algebraic subvariety $\mathcal{F}_{v} \subset U \times \GL_{m}$ by the equations
\[ Q(e_{j}, e_{\ell}) = \sum_{1 \leq i, k \leq m} a_{ij} a_{k\ell} \mathcal{Q}(v^{i}, v^{k}), \hspace{3em} 1 \leq j, \ell \leq m , \]
where $a_{ij}$ for $1 \leq i, j \leq m$ are the natural coordinates on $\GL_{m}$.

\begin{prop}
\label{alphares}
Denote by $\alpha : J^{d}_{r} U \times \GL_{m} \to J^{d}_{r} \ch{L}$ the map $J^{d}_{r} (q \circ \iota) \circ \beta$. Then if we restrict $\alpha$ to the fibre of $J^{d}_{r} U \times \GL_{m}$ above $\mathcal{F}_{v}$, the image of $\alpha$ lies in $J^{d}_{r} \ch{D}$. Moreover, if $\psi : B \to \ch{D}$ is a local period map with $B \subset U$ and decomposition $\psi = q \circ \iota \circ f$, and $\sigma \in J^{d}_{r} U$ with basepoint $s \in U$, we have $\alpha(\sigma, f(s)) = \psi \circ \sigma$.
\end{prop}

\begin{proof}
Suppose that $B \subset U$ is an analytic ball containing the point $s$ on which we have an analytic function $f : B \to \GL_{m}$ such that $b^{k} = \sum_{i = 1}^{m} f_{ik} v^{i}$ gives a flat basis. Then defining $\psi = q \circ \iota \circ f$ and applying \autoref{keymapprop2} we find that 
\[ \alpha(\sigma, f(s)) = q \circ \iota \circ \beta(\sigma, f(s)) = \psi \circ \sigma . \]
This in particular shows the second statement, so it suffices to show that the restriction of $\alpha$ to the fibre above $\mathcal{F}_{v}$ lands inside $J^{d}_{r} \ch{D}$. This is equivalent to arguing that if $f$ is chosen such that $(s, f(s))$ lies inside $\mathcal{F}_{v}$ for all $s \in B$, then $\psi(B)$ lies inside $\ch{D}$. Using the fact that $f^{-1}(s)$ gives the change-of-basis matrix from $v^{1}(s), \hdots, v^{m}(s)$ to $b^{1}(s), \hdots, b^{k}(s)$, we can observe that the condition that $(s, f(s))$ lies inside $\mathcal{F}_{v}$ amounts to the condition that $\mathcal{Q}(b^{i}(s), b^{j}(s)) = Q(e^{i}, e^{j})$. It then follows by \autoref{locpermapredef} that $\psi$ is in fact a local period map, hence $\alpha(\sigma, f(s)) = \psi \circ \sigma \in J^{d}_{r} \ch{D}$.
\end{proof}

\subsubsection{The maps $J^{d}_{r} S \to G \backslash J^{d}_{r} \ch{D}$}

We are now ready to construct the $K$-algebraic map $J^{d}_{r} S \to G \backslash J^{d}_{r} \ch{D}$ advertised in the introduction.

\jetcorresp*

\begin{proof}
Choose an open cover $\{ U_{i} \}_{i \in I}$ of $S$ by $K$-algebraic open sets on which $\mathcal{H}$ and the filtered pieces $F^{k} \mathcal{H}$ are trivial, and choose filtration-compatible frames $v^{1}_{i}, \hdots, v^{m}_{i}$ on each $U_{i}$. Associated to these frames we obtain maps $\alpha_{i} : J^{d}_{r} U_{i} \times \GL_{m} \to J^{d}_{r} \ch{L}$. Denote by $\pi_{i} : J^{d}_{r} U_{i} \times \GL_{m} \to U_{i} \times \GL_{m}$ the projection, and let $\alpha'_{i} : \pi_{i}^{-1}(\mathcal{F}_{v_{i}}) \to J^{d}_{r} \ch{D}$ be the restriction of $\alpha_{i}$ as in the statement of \autoref{alphares}. The $K$-algebraic variety $\pi^{-1}_{i}(\mathcal{F}_{v_{i}})$ is a $G$-torsor over $J^{d}_{r} U_{i}$, hence the $G$-invariant map $\alpha'_{i} : \pi^{-1}_{i}(\mathcal{F}_{v_{i}}) \to J^{d}_{r} \ch{D}$ induces a map of $K$-algebraic stacks $\eta^{d}_{r, i} : J^{d}_{r} U_{i} \to G \backslash J^{d}_{r} \ch{D}$. Let us note that if $j \in (J^{d}_{r} U_{i})(\mathbb{C})$ and $\psi : B \to \ch{D}$ is a local period map on $U_{i}$ at the base-point of $j$, then $\eta^{d}_{r,i}(j) = \psi \circ j$ modulo $G$ by \autoref{alphares}.

To complete the proof, it suffices to glue the maps $\eta^{d}_{r,i}$, which by the stacky formalism simply means gluing the $G$-torsors $q_{i} : \pi^{-1}_{i}(\mathcal{F}_{v_{i}}) \to J^{d}_{r} U_{i}$ to a $G$-torsor over $J^{d}_{r} S$, and gluing the maps $\alpha'_{i}$. If $v^{1}_{i}, \hdots, v^{m}_{i}$ and $v^{1}_{j}, \hdots, v^{m}_{j}$ are the two frames corresponding to the torsors $q_{i}$ and $q_{j}$, we have unique sections $b^{ij}_{kl}$ for $1 \leq k, l \leq m$ of the frame bundle associated to $\mathcal{H}$ such that $v^{k}_{i} = \sum_{m} b^{ij}_{kl} v^{l}_{j}$ on $U_{i} \cap U_{j}$. The uniqueness of the $b^{ij}_{kl}$ ensure that the induced maps between $q_{i}$ and $q_{j}$ on the overlap $U_{i} \cap U_{j}$ satisfy the cocycle condition, so we obtain a $G$-torsor over $J^{d}_{r} S$. The fact that the maps $\alpha'_{i}$ glue follows similarly.
\end{proof}

As a byproduct of \autoref{jetcorresp} we also obtain:

\begin{prop}
\label{torsorprop}
Let $p_{r} : \mathcal{T}^{d}_{r} \to J^{d}_{r} S$ be the $G$-torsor associated to the map $\eta^{d}_{r}$ of \autoref{jetcorresp}, and let $\alpha_{r} : \mathcal{T}^{d}_{r} \to J^{d}_{r} \ch{D}$ be the associated $G$-invariant map. Denote by $\pi^{r}_{0} : \mathcal{T}^{d}_{r} \to \mathcal{T}^{d}_{0}$ the natural projection. Then:
\begin{itemize}
\item[(i)] for each point $m_{0} \in \mathcal{T}^{d}_{0}$ there exists a local period map $\psi_{m_{0}} : B \to \ch{D}$ such that $p_{0}(m_{0}) \in B$;
\item[(ii)] if we regard local period maps up to equivalence of germs, then the action of $G$ on $\mathcal{T}_{0}$ is compatible with the action of $G$ on local period maps;
\item[(iii)] for any $j \in J^{d}_{r} S$ with basepoint $p_{0}(m_{0})$, there is a (necessarily unique) point $m_{r} \in \mathcal{T}^{d}_{r}$ projecting onto both $j$ and $m_{0}$, and such that $\alpha_{r}(m_{r}) = \psi_{m_{0}} \circ j$.
\end{itemize}
\end{prop}

\begin{proof}
All the properties are immediate from the construction of the torsor in the proof of \autoref{jetcorresp}, and may be checked locally on $S$. In particular, we may identify $\alpha_{r}$ with one of the maps $\pi_{i}^{-1}(\mathcal{F}_{v_{i}}) \to J^{d}_{r} \ch{D}$, the map $p_{r}$ with the projection $\pi^{-1}_{i}(\mathcal{F}_{v_{i}}) \to J^{d}_{r} U_{i}$, and $\mathcal{T}^{d}_{r}$ with $\pi_{i}^{-1}(\mathcal{F}_{v_{i}})$. Under these identifications, $m_{0}$ is identified with a point $(s, M_{0}) \in S(\mathbb{C}) \times \GL_{m}(\mathbb{C})$ that lies inside $\mathcal{F}_{v_{i}}$. For property (i) we construct $\psi_{m_{0}}$ as $q \circ \iota \circ f$ for $f$ solving the differential system associated to the frame $v_{i}$ with the initial condition $f(s) = M_{0}$. The action of $A \in G(\mathbb{C})$ on $\mathcal{T}^{d}_{0}$ then corresponds to the right-multiplication $M \mapsto M \cdot A$, which has the effect of replacing $\psi_{m_{0}}$ by $A^{-1} \cdot \psi_{m_{0}}$, which shows (ii). The point $m_{r}$ referred to in (iii) is the point $(j, M_{0})$, and so (iii) follows by \autoref{alphares}.
\end{proof}

\section{Hermitian Symmetric Foliations}
\label{hermfolsec}

We now turn our attention to the notion of \emph{Hermitian symmetric foliation} that is essential for \autoref{E3prop}(ii). Let us recall our setup. We fix an integral lattice $V$ and a polarizing form $Q : V \times V \to \mathbb{Z}$. We denote by $\ch{D}$ the flag variety parametrizing Hodge flags on $V$ satisfying the first Hodge-Riemann bilinear relation, $D \subset \ch{D}$ the open subset parametrizing Hodge flags additionally satisfying the positivity condition, and by $G = \textrm{Aut}(V, Q)$ the natural automorphism group. We suppose that we have a proper period map $\varphi : S \to \Gamma \backslash D$ with $\Gamma = G(\mathbb{Z})$, and we let $T = \varphi(S)$ be its image; this is a quasi-projective complex algebraic variety by \cite{OMINGAGA}. 

In what follows, we fix a lift $\widetilde{\varphi} : \widetilde{S} \to D \subset \ch{D}$ of $\varphi : S \to \Gamma \backslash D$, where $\pi : \widetilde{S} \to S$ is the universal cover, and denote by $\widetilde{T} = \widetilde{\varphi}(\widetilde{S})$ the image. We will denote by $\Gamma_{S} \subset \Gamma = G(\mathbb{Z})$ the subgroup obtained as the image of $\pi_{1}(S)$ under $\varphi_{*}$, and by $N \subset G$ the identity component of its Zariski closure. It follows from the structure theorem for period mappings in \cite[III.A]{GGK} that $\widetilde{T}$ lies inside a unique orbit $D_{T} \subset D$ of the identity component $N(\mathbb{R})^{\circ}$. We denote by $\ch{D}_{T} \subset \ch{D}$ the orbit of $N(\mathbb{C})$ containing $D_{T}$ as an open analytic subvariety. As we will be interested in discussing germs of $\widetilde{T}$, we note the following:

\begin{lem}
\label{closedT}
The locus $\widetilde{T}$ is a closed analytic subvariety of $D$, and hence of $D_{T}$.
\end{lem}

\begin{proof}
Choose a finite-index neat arithmetic subgroup $\Gamma' \subset \Gamma$. Then the action of $\Gamma'$ on $D$ is free and proper, and $\pi' : D \to \Gamma' \backslash D$ is a covering map. Choose a finite \'etale covering $g : S' \to S$ of $S$ such that $\varphi \circ g$ factors as $h \circ \varphi'$, where $\varphi' : S' \to \Gamma' \backslash D$ is a period map associated to $\mathbb{V}' = g^{*} \mathbb{V}$ and $h : \Gamma' \backslash D \to \Gamma \backslash D$ is the natural projection. As the map $\varphi \circ g$ is proper and $h$ is finite, it follows that $\varphi'$ is proper. We note that $\widetilde{\varphi}$ is also a lift of $\varphi'$.

Replacing $S$ with $S'$ we may then assume that the analytic projection $\pi : D \to \Gamma \backslash D$ is a covering map. Choose an open neighbourhood $U \subset D$ intersecting $\widetilde{T}$ such that $\restr{\pi}{U}$ is a homeomorphism. Then $\widetilde{T} \cap U$ is identified via $\pi$ with a component of $T \cap \pi(U)$, and is therefore closed analytic in $U$. Letting the open sets $U$ range over an open cover of $D$, we find that $\widetilde{T}$ is closed analytic in $D$.
\end{proof}

Let us now explain what we mean by \emph{Hermitian symmetric subvariety} and \emph{Hermitian symmetric foliation}.

\begin{defn}
\label{hermdef}
We say that a polarized variation of Hodge structure $\mathbb{V}$ on a smooth variety $S$ is \emph{Hermitian symmetric} if: 
\begin{itemize}
\item[(i)] in the case where the period map $\varphi$ is proper we have $\widetilde{T} = D_{T}$ in the above setup;
\item[(ii)] more generally for a possibly non-proper period map we have that $\widetilde{T}$ contains an open neighbourhood of $D_{T}$.
\end{itemize}
\end{defn}

\begin{rem}
The reader may wonder what the condition in \autoref{hermdef} has to do with Hermitian symmetric spaces. The point is that Hermitian symmetric domains are distinguished among period domains by the property that the Griffiths transversality condition is trivial; i.e., the Griffiths transverse subbundle is the entire tangent space. Since variations of Hodge structures satisfy Griffiths transversality, the conditions (i) and (ii) in \autoref{hermdef} necessarily imply that $D_{T}$ is a Hermitian symmetric domain, and allows us to think of $S$ as behaving similarly to a Shimura variety.
\end{rem}

\begin{defn}
We say that an irreducible complex subvariety $Y \subset S$ is Hermitian symmetric if the pullback $\restr{\mathbb{V}}{Y^{\textrm{sm}}}$ is Hermitian symmetric, where $Y^{\textrm{sm}} \to Y$ is a smooth resolution.\footnote{This is easily checked to be independent of the smooth resolution chosen.}
\end{defn}

\begin{defn}
\label{hermfoldef}
Given a variety $S$ with a polarized variation of Hodge structure $\mathbb{V}$, a \emph{Hermitian symmetric foliation} of $S$ is a dominant algebraic map $g : S \to U$ such that every irreducible component of every fibre $g^{-1}(u)$ is a Hermitian symmetric subvariety $Y \subset S$. It is called positive-dimensional if $Z = \varphi(Y)$ is positive-dimensional for each such $Y$.
\end{defn}

\subsection{Preliminaries from Functional Transcendence}

As explained in \autoref{folfromloc}, our argument for \autoref{E3prop}(ii) involves producing foliations of $\widetilde{T}$ by orbits of a Lie subgroup $H \subset G(\mathbb{C})$, and showing that these foliations are in fact the result of a Hermitian symmetric foliation of $S$ (or equivalently, of $T$). In fact, we will even show that the foliation induced by $H$ is algebraic, in the sense that it arises from intersecting $\widetilde{T}$ with an algebraic foliation of $\ch{D}_{T}$. We will then be left with the problem of  understanding the relationship between an algebraic foliation on $\ch{D}_{T}$ and an algebraic foliation on $S$, for which we will need the functional transcendence theory of variations of Hodge structures.

The crucial result is the following, due to Bakker and Tsimerman \cite{AXSCHAN}:

\begin{thm}[Ax-Schanuel]
\label{axschan}
Denote by $W \subset S \times \ch{D}_{T}$ the image of the graph of $\widetilde{\varphi}$, and let $V \subset S \times \ch{D}_{T}$ be an algebraic subvariety. Suppose that $U \subset V \cap W$ is an analytic component such that
\[ \textrm{codim}_{S \times \ch{D}_{T}} U < \textrm{codim}_{S \times \ch{D}_{T}} V + \textrm{codim}_{S \times \ch{D}_{T}} W . \]
Then the projection of $U$ to $S$ lies inside a proper weakly special subvariety.
\end{thm}
\begin{rem}
In \cite{AXSCHAN}, weakly special subvarieties are called \emph{weak Mumford-Tate}.
\end{rem}

\medskip
\medskip

Let us recall the notion of \emph{weakly special subvariety} that appears in the statement of \autoref{axschan}.

\begin{defn}
\label{weakspdef}
A weakly special subvariety of $S$ is an irreducible component of a subvariety of the form $\pi(\widetilde{\varphi}^{-1}(J \cdot t))$, where:
\begin{itemize}
\item[(i)] $t \in \widetilde{T}$ is a point;
\item[(ii)] $J$ is a $\mathbb{Q}$-normal subgroup of the Mumford-Tate group $\textrm{MT}(t)$ of $t$;
\item[(iii)] $\pi : \widetilde{S} \to S$ is the projection.
\end{itemize}
\end{defn}

\noindent For the fact that such subvarieties are algebraic, we refer to \cite{closurepositivelocus}. As is shown in \cite[Cor. 3.14]{closurepositivelocus}, weakly special subvarieties are precisely those irreducible subvarieties $Y \subset S$ that are maximal for their algebraic monodromy groups, i.e., those $Y$ for which the identity component of the Zariski closure of the monodromy representation on $\restr{\mathbb{V}}{Y^{\textrm{nor}}}$ is maximal, where $Y^{\textrm{nor}} \to Y$ is the normalization.

\medskip
\medskip

We record the consequences of the Ax-Schanuel Theorem which will important for us:

\begin{lem}
\label{zarclosurestruct}
The variety $\ch{D}_{T}$ is the Zariski closure of any analytic component of any germ of $\widetilde{T}$ in $\ch{D}$.
\end{lem}

\begin{proof}
As both $\ch{D}_{T}$ and $\widetilde{T}$ are irreducible, it suffices to show that $\widetilde{T}$ does not lie a subvariety of $\ch{D}$ of smaller dimension. If $\widetilde{T}$ lay in a proper subvariety $\ch{D}' \subset \ch{D}_{T}$ with $\dim \ch{D}' < \dim \ch{D}_{T}$, then the variety $V = S \times \ch{D}'$ would contain the image $W \subset S \times \ch{D}_{T}$ of the graph of $\widetilde{\varphi}$, hence we would have
\[ \codim_{S \times \ch{D}_{T}} W < \codim_{S \times \ch{D}_{T}} (S \times \ch{D}') + \codim_{S \times \ch{D}_{T}} W , \]
contradicting \autoref{axschan}.
\end{proof}

\begin{prop}
Suppose that $X \subset \ch{D}_{T}$ is an irreducible algebraic subvariety, that $t \in X \cap \widetilde{T}$ is a point, and we have a neighbourhood $B \subset D$ containing $t$ such that $B \cap X \subset B \cap \widetilde{T}$. Suppose $t = \widetilde{\varphi}(\widetilde{s})$, and let $C$ be an irreducible component of $\widetilde{\varphi}^{-1}(B \cap X)$ in $\widetilde{\varphi}^{-1}(B)$ whose image contains the germ $(X, t)$. Then $\pi(C)$ lies inside a Hermitian symmetric weakly special subvariety of $S$.
\label{liesinherm}
\end{prop}

\begin{proof}
We argue by induction on $\dim T$, the case where $\dim T = 0$ being automatic. The case where $\dim X = 0$ is also trivial, so we assume $\dim X > 0$. If $S$ is Hermitian symmetric there is nothing to prove, so we may assume $\dim \widetilde{T} < \dim \ch{D}_{T}$. We let $V = S \times X$.

Let $W \subset S \times \ch{D}_{T}$ be the projection of the graph of $\widetilde{\varphi}$. Let $s = \pi(\widetilde{s})$. Let $U$ be the analytic component of $W \cap V$ whose projection to $S$ contains $\pi(C, \widetilde{s})$. Then we have $(\dim U - \dim X) = (\dim C - \dim X) \geq (\dim S - \dim T)$, and hence
\begin{align*}
\textrm{codim}_{S \times \ch{D}_{T}} U &= \dim \ch{D}_{T} + \dim S - \dim U \\
&= (\dim \ch{D}_{T} - \dim X) + (\dim S - (\dim U - \dim X)) \\
&\leq \codim_{S \times \ch{D}_{T}} V + \dim \widetilde{T} \\
&< \codim_{S \times \ch{D}_{T}} V + \dim \ch{D}_{T} \\
&= \codim_{S \times \ch{D}_{T}} V + \codim_{S \times \ch{D}_{T}} W .
\end{align*}
By \autoref{axschan}, it follows that the projection of $U$ to $S$ lies a proper weakly special subvariety $Y \subset S$. By construction, it contains $\pi(C, \widetilde{s})$.

Let $Z = \varphi(Y)$ be the period image of $Y$, and let $Y^{\textrm{sm}} \to Y$ be a smooth resolution. Denote by $\varphi_{Y}$ the composition $Y^{\textrm{sm}} \to Y \xrightarrow{\varphi} Z \subset \Gamma \backslash D$. We may then construct an analogous diagram to the one associated to the variation on $S$:
\begin{center}
\begin{tikzcd}
\widetilde{Y^{\textrm{sm}}} \arrow[r, "\widetilde{\varphi}_{Y}"] \arrow[d, "\pi"] & \widetilde{Z} \arrow[d, "\pi"] \arrow[draw=none]{r}[sloped,auto=false,style={font=\normalsize}]{\subset} & \ch{D}_{Z} \\
Y^{\textrm{sm}} \arrow[r, "\varphi_{Y}"] & Z \arrow[draw=none]{r}[sloped,auto=false,style={font=\normalsize}]{\subset} & \Gamma \backslash D .
\end{tikzcd}
\end{center}
We may choose the lift $\widetilde{\varphi}_{Y}$ so that $\widetilde{Z}$ contains the point $t$, and hence so that $B \cap X \subset B \cap \widetilde{Z}$ for an appropriate neighbourhood $B \subset D$ (as a consequence of the fact that the projection of $(X, t)$ to $T$ lies inside $Z$). Replacing $C$ by an appropriate $C_{Y} \subset \widetilde{Y^{\textrm{sm}}}$ whose image contains $(X, t)$, we obtain the result from the induction hypothesis.
\end{proof}

\subsection{Preliminaries on $\Gamma_{S}$-orbits}

This section will be used to understand the implications of the Zariski density of $\Gamma_{S} \cap N(\mathbb{C})$ in $N(\mathbb{C})$ on analytic sets preserved by an analytic action of $\Gamma_{S} \cap N(\mathbb{C})$.

\begin{lem}
\label{serieslem}
Let $K_{0} = \Gamma_{S} \cap N(\mathbb{C})$, and denote by $K_{i} = [K_{i-1}, K_{i-1}]$ for $i \geq 1$. Then $\overline{K_{i}} = N$ for all $i$, where $\overline{(-)}$ denotes Zariski closure.
\end{lem}

\begin{proof}
Let $p: N_{1} \times \cdots \times N_{k} \to N$ be the $\mathbb{Q}$-simple decomposition of $N$. By \cite[III.A]{GGK}, we may assume the group $\Gamma_{S} \cap N(\mathbb{C})$ is the image under $p$ of an analogous decomposition $\Gamma_{1} \times \cdots \times \Gamma_{k}$ with $\overline{\Gamma_{i}} = N_{i}$. It therefore suffices to assume that $N$ is $\mathbb{Q}$-simple. We argue by induction on $i$, and let $L = \overline{K_{i}}$. Then as $K_{i}$ is stable under conjugation by $K_{i-1}$, so is $L$. As $\overline{K_{i-1}} = N$ by induction, the subgroup $L$ is stable under conjugation by a Zariski-dense subgroup of $N$, and therefore stable under conjugation by $N$ itself. It follows that $L$ is a normal non-trivial $\mathbb{Q}$-algebraic subgroup of $N$, hence $L = N$.
\end{proof}

\begin{lem}
\label{noanalyticsubgroup}
The group $\Gamma_{S} \cap N(\mathbb{C})$ does not lie in a proper analytic integral (in the sense of \cite[Def 1, pg.306]{intsub}) subgroup of $N(\mathbb{C})$.
\end{lem}

\begin{proof}
We define the $K_{i}$ as in \autoref{serieslem}, and suppose that $L_{0}$ is an analytic integral subgroup containing $K_{0} = \Gamma_{S} \cap N(\mathbb{C})$. If $L_{i}$ admits a non-trivial additive character $\tau_{i} : L_{i} \to \mathbb{C}$ then we let $L_{i+1}$ be its kernel, which is an analytic integral subgroup containing $K_{i+1}$. Eventually we must have that $L_{i}$ admits no non-trivial additive characters for some $i$, and that $\overline{L_{i}} = \overline{K_{i}} = N$ by \autoref{serieslem}. It then follows from \cite[Prop. 5]{magid} that $L_{i} = N$.
\end{proof}

\subsection{Germs of Local Symmetries}
\label{locsymgermsec}

In this section we study the sets $H_{t}$, defined by:

\begin{defn}
\label{Htdef}
For $t \in \widetilde{T}$, we write $H_{t} \subset G(\mathbb{C})$ for the set
\[ H_{t} = \{ g \in G(\mathbb{C}) : g \cdot (\widetilde{T}, t) \subset (\widetilde{T}, g \cdot t) \} . \]
\end{defn}

\noindent As $t \in \widetilde{T}$ varies, the sets $H_{t}$ contain all those elements $g \in G(\mathbb{C})$ which are responsible for equating points of $S(\mathbb{C})$ under the equivalence relation (E3). We will argue that the germ at the identity of $H_{t}$ is independent of $t$, and coincides with the germ of a connected Lie subgroup $H \subset G(\mathbb{C})$ which induces a foliation of $\widetilde{T}$. We will later relate the condition that the map of \autoref{E3prop}(ii) be quasi-finite to the condition that $H$ be the trivial group.

\medskip

\begin{lem}
\label{basicHfacts}
We have
\begin{itemize}
\item[(i)] if $g \in H_{t}$, then $g \cdot \ch{D}_{T} = \ch{D}_{T}$;
\item[(ii)] if $\gamma \in \Gamma_{S}$, then $\gamma \cdot H_{t} = H_{t}$;
\item[(iii)] if $\gamma \in \Gamma_{S}$, then $H_{t} \cdot \gamma = H_{\gamma^{-1} \cdot t}$. 
\end{itemize}
\end{lem}

\begin{proof}
For (i) we use that by \autoref{zarclosurestruct} above, the Zariski closures of any two components of the germs $(\widetilde{T}, t)$ and $(\widetilde{T}, g \cdot t)$ agree and equal $\ch{D}_{T}$, and the fact that $g$ acts by an algebraic automorphism. Property (ii) and (iii) are immediate from the definition and the fact that $\Gamma_{S}$ preserves $\widetilde{T}$.
\end{proof}

\medskip

The next three Lemmas show that $H_{t}$ is independent of $t$ near the identity.

\begin{lem}
\label{Hext}
Let $t \in \widetilde{T}$, and $\alpha \in H_{t}$. Then there exists neighbourhoods $B \subset D$ of $t$ and $\mathcal{G} \subset G(\mathbb{C})$ of $\alpha$ such that 
\[ \mathcal{G} \cap H_{t} \subset \mathcal{G} \cap H_{t'} , \]
for any $t' \in \widetilde{T} \cap B$. If in addition $(\widetilde{T}, t)$ an irreducible analytic germ, equality holds for all such $t'$.
\end{lem}

\begin{proof}
Fix a ball $B \subset D$ around $t \in \widetilde{T}$ such that the irreducible components of $(\widetilde{T}, t)$ extend to irreducible analytic sets $C_{1}, \hdots, C_{r}$ in $B$ giving a decomposition $\widetilde{T} \cap B = C_{1} \cup \cdots \cup C_{r}$. Shrinking $B$ if necessary, let $\mathcal{G}$ be a neighbourhood of $\alpha$ such that $g \cdot B \subset D$ for all $g \in \mathcal{G}$.\footnote{One can find such a neighbourhood by choosing $B$ to have compact closure and using the fact that the topology on the automorphism group of $\textrm{Aut}(\ch{D})$ coincides with the compact-open topology.} Then as $g \cdot C_{i}$ is irreducible in $g \cdot B$ for each $g \in \mathcal{G}$ and each $i$, and $\widetilde{T} \cap (g \cdot B)$ is closed in $g \cdot B$ by our choice of $B$, we learn that $g \cdot C_{i} \subset \widetilde{T}$ as soon as $g \cdot (C_{i}, t') \subset (\widetilde{T}, g \cdot t')$ for some $t' \in C_{i}$. In particular, we obtain $\mathcal{G} \cap H_{t} \subset \mathcal{G} \cap H_{t'}$ for any $t' \in \widetilde{T} \cap B$. If $(\widetilde{T}, t)$ is irreducible then $r = 1$ and $C = C_{1} = \widetilde{T} \cap B$, meaning that $g \cdot (\widetilde{T}, t') \subset (\widetilde{T}, g \cdot t')$ necessarily implies that $g \cdot (\widetilde{T}, t) \subset (\widetilde{T}, g \cdot t)$ as well.
\end{proof}

\begin{lem}
\label{Htclosedanal}
Let $U \subset \ch{D}$ be any open subset such that $\widetilde{T} \cap U$ is a closed analytic subset of $U$. Denote by $\mathcal{G} \subset G(\mathbb{C})$ the open subset of $g \in G(\mathbb{C})$ such that $g \cdot t \in U$. Then $H_{t} \cap \mathcal{G}$ is a closed analytic subset of $\mathcal{G}$.
\end{lem}

\begin{proof}
Fix an analytic component $(C, t) \subset (\widetilde{T}, t)$. It suffices to show that the condition on $g \in G(\mathbb{C})$ defined by $g \cdot (C, t) \subset (\widetilde{T}, g \cdot t)$ is an analytic condition on $\mathcal{G}$. For this we may fix a compatible family (recall \autoref{compjetdef}) of non-degenerate jets $j_{r} \in J^{d}_{r,nd} C$ with $j_{0} = t$ and $d = \dim C$. Applying \autoref{jetdimlem} to $(g \cdot C) \cap \widetilde{T}$ we find that the condition that $g \cdot (C, t) \subset (\widetilde{T}, g \cdot t)$ is equivalent to the condition that $g \cdot j_{r} \in J^{d}_{r} \widetilde{T}$ for all $r$. As the map $\mathcal{G} \times \{ j_{r} \} \to J^{d}_{r} U$ given by $g \mapsto g \cdot j_{r}$ is analytic, this condition is given by analytic conditions on $\mathcal{G}$.
\end{proof}

\begin{notn}
We denote by $\widetilde{T}_{\textrm{irr}}$ the locus of $t \in \widetilde{T}$ such that $(\widetilde{T}, t)$ is irreducible as an analytic germ.
\end{notn}

\begin{lem}
\label{Hindep}
The germ of $H_{t} \subset G(\mathbb{C})$ at the identity is independent of $t$.
\end{lem}

\begin{proof}
As $\widetilde{T}_{\textrm{irr}}$ is an open analytic subset of $\widetilde{T}$, it is necessarily connected since $\widetilde{T}$ is the image of the irreducible analytic variety $\widetilde{S}$. Fix $t_{0} \in \widetilde{T}_{\textrm{irr}}$, and consider the condition on $t \in \widetilde{T}_{\textrm{irr}}$ given by the equality of germs $(H_{t}, \textrm{id}) = (H_{t_{0}}, \textrm{id})$. Taking $\alpha = \textrm{id}$ \autoref{Hext} above shows that this is an open condition on $t \in \widetilde{T}_{\textrm{irr}}$; we claim it is also closed. Indeed, supposing we have a convergent sequence $t_{i} \to t$ with $t_{i}, t \in \widetilde{T}_{\textrm{irr}}$ and $(H_{t_{i}}, \textrm{id}) = (H_{t_{0}}, \textrm{id})$, it suffices to show that $(H_{t}, \textrm{id}) = (H_{t_{i}}, \textrm{id})$ for some $i$. But applying \autoref{Hext} again the germ $(H_{t}, \textrm{id})$ is constant in a neighbourhood of $t$, hence must equal the germ $(H_{t_{i}}, \textrm{id})$ for large enough $i$.

Considering now the case of a general $t \in \widetilde{T}$, it suffices by \autoref{Hext} to show that $(H_{t'}, \textrm{id}) \subset (H_{t}, \textrm{id})$ for some $t' \in \widetilde{T}_{\textrm{irr}}$. Let us enumerate the components $(C_{1}, t), \hdots, (C_{r}, t)$ of $(\widetilde{T}, t)$, and extend them to components $C_{1}, \hdots, C_{r}$ in a neighbourhood $B$ of $t$. We then choose $t'_{i} \in C_{i} \cap \widetilde{T}_{\textrm{irr}}$ for each $i$. By the above reasoning $(H_{t'_{1}}, \textrm{id}) = \cdots = (H_{t'_{r}}, \textrm{id})$. Hence if $\mathcal{G} \subset G(\mathbb{C})$ is a sufficiently small neighbourhood of the identity and $g \in \mathcal{G} \cap H_{t'_{1}}$ we find that $g \cdot C_{i} \subset \widetilde{T}$ for each $i$, which implies that $g \in H_{t}$.
\end{proof}

\begin{cor}
\label{conjstable}
The germ $(H_{t}, \textrm{id})$ is invariant under conjugation by $\Gamma_{S}$.
\end{cor}

\begin{proof}
Choose $\gamma \in \Gamma_{S}$. By \autoref{basicHfacts} above we have
\[ \gamma \cdot (H_{t}, \textrm{id}) \cdot \gamma^{-1} = (H_{\gamma^{-1} \cdot t}, \textrm{id}) , \]
but this is equal to $(H_{t}, \textrm{id})$ by \autoref{Hindep}.
\end{proof}

\medskip

\medskip

We now turn to the issue of showing that the germs $(H_{t}, \textrm{id})$ agree with that of a connected integral Lie subgroup of $G(\mathbb{C})$.

\begin{lem}
\label{locmultlem}
For any positive integer $m$, there exists a neighbourhood $U$ of the identity in $G(\mathbb{C})$ such that if $h_{1}, \hdots, h_{m} \in U \cap H_{t}$ then $h_{1} h_{2} \cdots h_{m} \in H_{t}$.
\end{lem}

\begin{proof}
It suffices to show this for $m = 2$. Fix a ball $B$ and $\mathcal{G}$ as in \autoref{Hext} with $\alpha = \textrm{id}$. We may choose $U$ to be a neighbourhood of the identity contained in $\mathcal{G}$ such that if $u \in U$ then $u \cdot t \in B$. Then if $h_{1}, h_{2} \in U \cap H_{t}$ we will have that $h_{1} \in H_{h_{2} \cdot t}$, hence $h_{1} h_{2} \in H_{t}$.
\end{proof}

\begin{prop}
\label{isaLiegroup}
The germ $(H_{t}, \textrm{id})$ is the germ of a connected integral analytic Lie subgroup of $G(\mathbb{C})$.
\end{prop}

\begin{proof}
Denote by $\mathfrak{g}_{\mathbb{C}}$ the Lie algebra of $G(\mathbb{C})$, and let $\mathfrak{h} \subset \mathfrak{g}_{\mathbb{C}}$ be the complex subspace corresponding to the tangent space of $H_{t}$. Denote by $\textrm{pow}_{k} : G(\mathbb{C}) \to G(\mathbb{C})$ the map $g \mapsto g^{k}$. It follows from \autoref{locmultlem} above that $\textrm{pow}_{k}$ preserves the germ $(H_{t}, \textrm{id})$ for all $k$. Let $\exp : \mathfrak{g}_{\mathbb{C}} \to G(\mathbb{C})$ be the usual exponential map, and choose a logarithm $\log : U \to \mathfrak{g}_{\mathbb{C}}$ in a neighbourhood $U$ of the origin. Translating the above we find that the germ $(\log(H_{t}), 0)$ is stable under scalar multiplication by $\mathbb{Z}$. 

We claim that any analytic germ $(A, 0)$ at the origin of a complex vector space $V$ and which is stable under $\mathbb{Z}$-multiplication must be a linear subspace; if so, then $(\log(H_{t}), 0) = (\mathfrak{h}, 0)$ for some linear subspace $\mathfrak{h}$, completing the proof. To see this, identify $V = \mathbb{C}^{m}$ and let $f(z_{1}, \hdots, z_{m})$ be a power series centred at $0$ vanishing on $(A, 0)$. Then by assumption, $f(k z_{1}, \hdots, k z_{m})$ also vanishes on $(A, 0)$ for all $k \in \mathbb{Z}$. By taking linear combinations of $f(k z_{1}, \hdots, k z_{m})$ for $k \in \mathbb{Z}$ we may inductively construct analytic functions $f_{i}$ with $f_{1} = f$ and such that $f_{i}$ has the same linear term as $f_{i-1}$ but no non-linear terms before degree $i+1$; so for instance, $f_{2}$ has no quadratic terms and may be defined by 
\[ f_{2}(z_{1}, \hdots, z_{m}) = \frac{1}{2} \left(4 f(z_{1}, \hdots, z_{m}) - f(2 z_{1}, \hdots, 2 z_{m}) \right) . \]
Denoting by $L$ the linear part of $f$ and $W$ the subspace defined by $L = 0$, it then follows that the fibre $(J^{d}_{r} A)_{0}$ lies inside $(J^{d}_{r} W)_{0}$ for all $d$ and $r$; we therefore have that $(A, 0) \subset (W, 0)$ by \autoref{jetdimlem} applied to $(A \cap W, 0)$. We then replace $V$ with $W$ and argue by induction.
\end{proof}

\medskip

As a consequence of the preceding results, we may now define:

\begin{defn}
\label{Hdef}
We denote by $H \subset G(\mathbb{C})$ the connected integral Lie subgroup whose tangent space at the identity agrees with the germ $(H_{t}, \textrm{id})$ for all $t \in \widetilde{T}$.
\end{defn}

\begin{lem}
\label{stableunderN}
The group $H$ is stable under conjugation by $N(\mathbb{C})$. 
\end{lem}

\begin{proof}
The subgroup of $N(\mathbb{C})$ defined by the condition $n \cdot H \cdot n^{-1} = H$ is an analytic integral subgroup which contains $\Gamma_{S} \cap N(\mathbb{C})$ by \autoref{conjstable}. By \autoref{noanalyticsubgroup} it must equal all of $N(\mathbb{C})$.
\end{proof}

\subsection{The Elements Inducing (E3)}

In this section we give a characterization of the equivalence relation (E3) in terms of the sets $H_{t}$; in particular, we show that the elements $g \in G(\mathbb{C})$ which equate points of $S(\mathbb{C})$ under (E3) may be taken to lie in the union $\bigcup_{t \in \widetilde{T}} H_{t}$, and that this union lies inside a countable set of cosets of $H$. These results will not be needed until \autoref{equivrelsec} and can be temporarily skipped.

\begin{lem}
\label{Htlocalstructure}
For any $t \in \widetilde{T}$ and $g \in H_{t}$, the analytic germ $(H_{t}, g)$ is irreducible of dimension equal to $d = \dim H$.
\end{lem}

\begin{proof}
We first observe that 
\[ (H, \textrm{id}) \cdot g = (H_{g \cdot t}, \textrm{id}) \cdot g \subset (H_{t}, g) , \]
hence $\dim (H_{t}, g) \geq d$. It therefore suffices to show the irreducibility of $(H_{t}, g)$.

We claim that we may reduce to the case where $g \cdot t \in \widetilde{T}_{\textrm{irr}}$. If we apply \autoref{Hext} to $\alpha = g$, we obtain neighbourhoods $B$ of $t$ and $\mathcal{G}$ of $g$ such that
\[ \mathcal{G} \cap H_{t} \subset \mathcal{G} \cap H_{t'} , \]
for any $t' \in \widetilde{T} \cap B$, hence in particular $(H_{t}, g) \subset (H_{t'}, g)$ for such $t'$. As $(H_{t'}, g)$ being irreducible of dimension $d$ combined with $\dim (H_{t}, g) \geq d$ will imply the same fact for $(H_{t}, g)$, we may use the density of $\widetilde{T}_{\textrm{irr}}$ in $\widetilde{T}$ to replace $t$ with $t' \in (g^{-1} \cdot \widetilde{T}_{\textrm{irr}}) \cap B$.

From the fact that $(\widetilde{T}, g \cdot t)$ is irreducible we learn that $g \cdot (\widetilde{T}, t) = (\widetilde{T}, g \cdot t)$ (instead of just possibly the inclusion $g \cdot (\widetilde{T}, t) \subset (\widetilde{T}, g \cdot t)$). Thus we have that $g^{-1} \in H_{g \cdot t}$. We then find that
\[ (H_{t}, g) \cdot g^{-1} \subset (H_{g \cdot t}, \textrm{id}) = (H, \textrm{id}) , \]
so the result follows by the irreducibility of the germ $(H, \textrm{id})$.
\end{proof}

\begin{cor}
\label{irredcor}
For each $t \in \widetilde{T}$ and $g \in H_{t}$, we have $(H_{t}, g) = (H \cdot g, g)$.
\end{cor}

\begin{proof}
We have an inclusion $(H_{g \cdot t}, \textrm{id}) \cdot g \subset (H_{t}, g)$, so the result follows using the equality $(H_{g \cdot t}, \textrm{id}) = (H, \textrm{id})$ from \autoref{Hindep} and the result of \autoref{Htlocalstructure} above.
\end{proof}

\begin{notn}
We define
\[ \widetilde{H} = \bigcup_{t \in \widetilde{T}} H_{t} . \]
\end{notn}

\begin{lem}
\label{countrightcoset}
The set $\widetilde{H} \subset G(\mathbb{C})$ lies inside a countable union of right $H$-cosets.
\end{lem}

\begin{proof}
Suppose we have $g \in \widetilde{H}$. Then $g$ lies in $H_{t}$ for some $t$. Applying \autoref{Hext} we obtain neighbourhoods $B_{t, g}$ of $t$ and $\mathcal{G}_{t, g}$ of $g$ such that
\[ H_{t} \cap \mathcal{G}_{t, g} = H_{t'} \cap \mathcal{G}_{t, g} = (H \cdot g) \cap \mathcal{G}_{t, g} \]
for any $t' \in B_{t, g}$, where we use the local description of the germs of $H_{t'}$ afforded by \autoref{irredcor} to obtain the equalities. The union $\mathcal{U} = \bigcup_{t, g} B_{t, g} \times \mathcal{G}_{t, g}$ is open in $\ch{D} \times G(\mathbb{C})$, hence by the second countability of the topology we can find a countable subcover $B_{t_{i}, g_{i}} \times \mathcal{G}_{t_{i}, g_{i}}$ so that $\mathcal{U} = \bigcup_{i = 1}^{\infty} B_{t_{i},g_{i}} \times \mathcal{G}_{t_{i}, g_{i}}$. We now claim that
\begin{equation}
\label{countableunion}
\widetilde{H} \subset \bigcup_{i = 1}^{\infty} (H \cdot g_{i}) \cap \mathcal{G}_{t_{i}, g_{i}}
\end{equation}
which will give the result. Indeed, given $g' \in \widetilde{H}$ we have $g' \in H_{t'}$ for some $t'$. Then $(t', g') \in \mathcal{U}$, so $(t', g') \in B_{t_{i}, g_{i}} \times \mathcal{G}_{t_{i}, g_{i}}$ for some $i$. This means that 
\[ H_{t'} \cap \mathcal{G}_{t_{i}, g_{i}} = (H \cdot g_{i}) \cap \mathcal{G}_{t_{i}, g_{i}} , \]
hence $g'$ lies in the right-hand side of \autoref{countableunion}. 
\end{proof}

\begin{prop}
\label{E3redef}
There exists a countable sequence $\{ g_{i} \}_{i = 1}^{\infty}$ with $g_{i} \in G(\mathbb{C})$ such that the relation (E3) on $S(\mathbb{C})$ in \autoref{equivdef} is equivalent to the following:
\begin{itemize}
\item[(E3)'] there exists distinguished local lifts $\psi_{i} : B_{i} \to \widetilde{T} \subset D$ with $s_{i} \in B_{i}$ for $i = 1, 2$ and $g \in \bigcup_{i = 1}^{\infty} H \cdot g_{i}$ such that $g \cdot \mathcal{C}(\psi_{1}, s_{1}) = \mathcal{C}(\psi_{2}, s_{2})$.
\end{itemize}
\end{prop}

\begin{proof}
It follows immediately from the definition of (E3) that we may restrict to just those local lifts landing in $\widetilde{T}$. Then any $g \in G(\mathbb{C})$ such that $g \cdot \mathcal{C}(\psi_{1}, s_{1}) = \mathcal{C}(\psi_{2}, s_{2})$ will necessarily send some irreducible germ of $\widetilde{T}$ at a point $t$ near $\psi_{1}(s_{1})$ to an irreducible germ of $\widetilde{T}$ near $\psi_{2}(s_{2})$, hence lie inside $H_{t}$ for some $t$. Thus we obtain the same equivalence relation considering only those $g \in G(\mathbb{C})$ that lie inside $\widetilde{H}$, and hence the equivalence relation is also unchanged by considering those $g$ lying inside the countable union $\bigcup_{i = 1}^{\infty} H \cdot g_{i}$ containing $\widetilde{H}$ that exists by \autoref{countrightcoset}.
\end{proof}

\subsection{The Foliation Determined by $H$}

The preceding sections establish that the germs of $\widetilde{T}$ are locally foliated by the orbits of the group $H$. In this section we employ tools from the functional transcendence theory of variations of Hodge structures to establish the global algebraic properties of this foliation. 

In what follows we denote by $M$ the identity component  $\textrm{Aut}(\ch{D}_{T})^{\circ}$ of the algebraic group $\textrm{Aut}(\ch{D}_{T})$. By \autoref{basicHfacts}(i) above we know that $H$ acts on $\ch{D}_{T}$, hence acts through a subgroup of $M$. We hence denote by $i(N)$ (resp. $i(H)$) the image of $N$ in $M$ (resp. the image of $H$ in $M(\mathbb{C})$). We let $i(N) = N_{1} \times \cdots \times N_{k}$ and $M = M_{1} \times \cdots \times  M_{\ell}$ be the decompositions into simple factors; we note that the theory of flag varieties (c.f. \cite[Ch. 3.3]{flagvar}) ensures that we may choose $k = \ell$ and $N_{i} \subset M_{i}$ for each $i$. We analogously write  $i(\mathfrak{n}_{\mathbb{C}}) = \mathfrak{n}_{1} \oplus \cdots \oplus \mathfrak{n}_{k}$ and $\mathfrak{m}_{\mathbb{C}} = \mathfrak{m}_{1} \oplus \cdots \oplus \mathfrak{m}_{k}$ for the associated Lie algebra decompositions.

\begin{lem}
\label{Hisprod}
After possibly replacing the almost-direct decompositions of $N$ and $M$ with equivalent ones, the group $i(H)$ is the set of complex points of a product
\[ H_{1} \times \cdots \times H_{k} , \]
of simple complex algebraic groups, where for each $i$, $H_{i}$ belongs to the set $\{ \{ 1 \}, N_{i}, M_{i} \}$.
\end{lem}

\begin{proof}
The theory of complex flag varieties \cite[Ch. 3.3]{flagvar} ensures that we have $\mathfrak{n}_{i} \subset \mathfrak{m}_{i}$, with equality in all but three possible exceptional cases:
\begin{itemize}
\item[(i)] the inclusion $\mathfrak{n}_{i} \hookrightarrow \mathfrak{m}_{i}$ is isomorphic to $\mathfrak{sp}_{2n} \hookrightarrow \mathfrak{sl}_{2n}$ for some $n \geq 2$;
\item[(ii)] the inclusion $\mathfrak{n}_{i} \hookrightarrow \mathfrak{m}_{i}$ is isomorphic to $\mathfrak{so}_{2n+1} \hookrightarrow \mathfrak{so}_{2n+2}$ for some $n \geq 3$;
\item[(ii)] the inclusion $\mathfrak{n}_{i} \hookrightarrow \mathfrak{m}_{i}$ is isomorphic to $\mathfrak{g}_{2} \hookrightarrow \mathfrak{so}_{7}$.
\end{itemize}
In each case, $\mathfrak{n}_{i}$ is a maximal proper Lie subalgebra of $\mathfrak{m}_{i}$. By \autoref{stableunderN} above, the set $L = i(N) \cdot i(H)$ is a Lie subgroup of $M(\mathbb{C})$ with Lie algebra $\mathfrak{l}$ satisfying $i(\mathfrak{n}_{\mathbb{C}}) \subset \mathfrak{l} \subset \mathfrak{m}$. The ideal $\mathfrak{l}_{i}$ of $\mathfrak{l}$ generated by $\mathfrak{n}_{i} \subset \mathfrak{l}$ is a subalgebra of $\mathfrak{m}_{i}$ containing $\mathfrak{n}_{i}$, hence $\mathfrak{l}_{i} \in \{ \mathfrak{n}_{i}, \mathfrak{m}_{i}\}$. The projections $\mathfrak{l} \to \mathfrak{m}_{i}$ have image equal to $\mathfrak{l}_{i}$, hence $\mathfrak{l}$ is semisimple with decomposition $\mathfrak{l} = \mathfrak{l}_{1} \oplus \cdots \oplus \mathfrak{l}_{k}$. The result now follows from the fact that $i(H)$ is normal in $L$ by \autoref{stableunderN}, hence the Lie algebra $i(\mathfrak{h})$ of $i(H)$ is a sum of simple summands of $\mathfrak{l}$.
\end{proof}

\begin{cor}
\label{alggroupcor}
The group $H \subset G(\mathbb{C})$ is a complex algebraic subgroup.
\end{cor}

\begin{proof}
The subgroup $G_{T} \subset G$ of automorphisms preserving $\ch{D}_{T}$ is algebraic, and \autoref{Hisprod} shows that $H$ may be identified with the inverse image of $H_{1} \times \cdots \times H_{k}$ under the algebraic map $G_{T} \to M$.
\end{proof}

\begin{lem}
\label{equivQgroup}
There exists a normal $\mathbb{Q}$-algebraic subgroup $J \subset N$ such that for each $t \in \widetilde{T}$ we have an equality of orbits $H \cdot t = J \cdot t$.
\end{lem}

\begin{proof}
We begin by observing that it suffices to show that there exists $J$ such that $H \cdot t = J \cdot t$ for one such $t$. Indeed by the theory of flag varieties we have $\ch{D}_{T} \simeq \ch{D}_{1} \times \cdots \times \ch{D}_{k}$ in accordance with the simple factors $M_{1}, \hdots, M_{k}$ of $\textrm{Aut}(\ch{D}_{T})^{\circ}$. As the orbits of $H$ and $J$ correspond to orbits of products of the $M_{i}$, the equality $H \cdot t_{0} = J \cdot t_{0}$ for one $t_{0} \in \widetilde{T}$ necessarily implies that $H$ and $J$ act through the same set of factors, so the result follows.

Now let us unravel the consequences of \autoref{liesinherm}. Given a point $t \in \widetilde{T}$, we denote by $\textrm{MT}(t)$ its Mumford-Tate group. We let $\textrm{MT}(S)$ denote the Mumford-Tate group at a generic point of $\widetilde{T}$. We define $\mathcal{W} \subset \widetilde{T}$ to be the countable union of analytic subvarieties of $\widetilde{T}$ of the following two types:
\begin{itemize}
\item[(i)] components of $\widetilde{T} \cap (\textrm{MT}(t) \cdot t)$ where the containment $\textrm{MT}(t) \subsetneq \textrm{MT}(S)$ is proper;
\item[(ii)] the analytic loci 
\[ \mathcal{A}_{J} = \{ t \in \widetilde{T} : \dim_{t} (\widetilde{T} \cap (J \cdot t)) = \dim (J \cdot t) \} , \]
for $J \subset N$ a normal $\mathbb{Q}$-algebraic subgroup such that $H \cdot t_{0} \subset J \cdot t_{0}$ for some (hence any) $t_{0} \in \widetilde{T}$.
\end{itemize}
We claim that $\mathcal{W} = \widetilde{T}$. If not, we may take a point $t_{0} \in \widetilde{T} \setminus \mathcal{W}$, and apply \autoref{liesinherm} to the germ $(H \cdot t_{0}, t_{0}) \subset (\widetilde{T}, t_{0})$. We learn that there is a Hermitian symmetric weakly special subvariety $Y \subset S$ such that $(\pi(H \cdot t_{0}), \pi(t_{0})) \subset (\varphi(Y), \pi(t_{0}))$. Choose a smooth resolution $Y^{\textrm{sm}} \to Y$ and construct the diagram 
\begin{equation}
\label{Ylift}
\begin{tikzcd}
\widetilde{Y^{\textrm{sm}}} \arrow[r, "\widetilde{\varphi}_{Y}"] \arrow[d, "\pi"] & \widetilde{Z} \arrow[d, "\pi"] \arrow[draw=none]{r}[sloped,auto=false,style={font=\normalsize}]{\subset} & \ch{D}_{Z} \\
Y^{\textrm{sm}} \arrow[r, "\varphi_{Y}"] & Z \arrow[draw=none]{r}[sloped,auto=false,style={font=\normalsize}]{\subset} & \Gamma \backslash D ,
\end{tikzcd}
\end{equation}
where $\widetilde{Z}$ is chosen to pass through $t$. The fact that $(\pi(H \cdot t_{0}), \pi(t_{0})) \subset (\varphi(Y), \pi(t_{0}))$ then guarantees that $(H \cdot t_{0}, t_{0}) \subset (\widetilde{Z}, t_{0})$. As $t_{0}$ was chosen so that $\textrm{MT}(t_{0}) = \textrm{MT}(S)$, the algebraic monodromy group $N_{Y}$ generating $\ch{D}_{Z}$ must be normal in $N$. As $Y$ is Hermitian symmetric we must have $\dim_{t_{0}} \widetilde{Z} = \dim_{t_{0}} \ch{D}_{Z} = \dim_{t_{0}} (N_{Y} \cdot t_{0})$, contradicting the claim that $t_{0}$ does not lie in $\mathcal{W}$.

It follows from the irreducibility of $\widetilde{T}$ that one of the constituent analytic sets that make up $\mathcal{W}$ must equal all of $\widetilde{T}$; the sets of type (i) are all proper subsets, hence we have $\widetilde{T} = \mathcal{A}_{J}$ for a $\mathbb{Q}$-algebraic normal subgroup $J \subset N$. Letting $D_{1} \subset D_{T}$ be the orbit of $J(\mathbb{R})^{\circ}$, we obtain as in \cite[III.A]{GGK} a factorization $D_{T} = D_{1} \times D_{2}$. Choose a product of open balls $B_{1} \times B_{2}$ such that $\widetilde{T} \cap (B_{1} \times B_{2})$ is an irreducible closed analytic subset. From irreducibility and the fact that $\mathcal{A}_{J} = \widetilde{T}$ we learn that $B_{1} \times \{ t_{2} \} \subset \widetilde{T}$ for any point $(t_{1}, t_{2}) \in \widetilde{T}$. Denote by $\widetilde{T}_{2}$ the projection of $\widetilde{T} \cap (B_{1} \times B_{2})$ to $B_{2}$. If we choose $g \in J(\mathbb{C})$ sufficiently close to the identity and $t \in \widetilde{T} \cap (B_{1} \times B_{2})$ we will then have
\[ g \cdot (\widetilde{T} \cap (B_{1} \times B_{2}), t) = (g \cdot (B_{1} \times \widetilde{T}_{2}), g \cdot t) = ((g \cdot B_{1}) \times \widetilde{T}_{2}, g \cdot t) = (\widetilde{T} \cap (B_{1} \times B_{2}), g \cdot t) . \]
It follows that the germ of $J$ at the identity lies in $H_{t}$, hence $J \subset H$.
\end{proof}

\begin{cor}
\label{hermsym}
For $t \in \widetilde{T}$ let $\mathcal{W}_{t} = (H \cdot t) \cap \widetilde{T} = (J \cdot t) \cap \widetilde{T}$ (where $J$ is as in \autoref{equivQgroup}) be the associated orbit. Then $\pi(\widetilde{\varphi}^{-1}(\mathcal{W}_{t}))$ is a union of Hermitian symmetric weakly special subvarieties of $S$.
\end{cor}

\begin{proof}
The fact that these varieties are weakly special is true by definition, so let us fix such a variety $Y \subset S$ and show it is Hermitian symmetric. As the orbits of $H$ and $J$ agree we have $(J \cdot t, t) \subset (\widetilde{T}, t)$ for each $t \in \widetilde{T}$. Constructing the diagram (\ref{Ylift}) above for this choice of $Y$ we find that $\dim_{t} \widetilde{Z} = \dim J \cdot t$. By \autoref{zarclosurestruct} we must have $J \cdot t = \ch{D}_{Z}$, hence the result.
\end{proof}

\medskip

We finally conclude with the main result of this section:

\begin{defn}
Given a variation of Hodge structure $\mathbb{V}$ on $S$ and $N \subset \textrm{MT}(S)$ a normal subgroup of the Mumford-Tate group $\textrm{MT}(S)$, we denote by $\mathbb{V}^{N}$ the variation of $N$-invariants. 
\end{defn}

\begin{prop}
\label{hermfolprop}
Let $M = \textrm{MT}(S)$ be the Mumford-Tate group of $\mathbb{V}$, and let $J_{1} = J$ be the subgroup of $N$ given by  \autoref{equivQgroup}. Let $J_{2} \subset N$ and $N' \subset M$ be complementary normal factors giving almost-direct decompositions $N = J_{1} \cdot J_{2}$ and $M = N \cdot N'$. Then there exists finitely many integers $a_{ij}, b_{ij}$ with $1 \leq i \leq m$ and $j = 1, 2$ such that the variations of invariants
\begin{align*}
\mathbb{V}_{j} := \left(\bigoplus_{i} \mathbb{V}^{\otimes a_{ij}} \otimes (\mathbb{V}^{\vee})^{\otimes b_{ij}}\right)^{J_{j} \cdot N'} ;
\end{align*}
have the following properties:
\begin{itemize}
\item[(i)] Denote by $V_{j} = \left(\bigoplus_{i} V^{\otimes a_{ij}} \times (V^{\vee})^{\otimes b_{ij}}\right)^{J_{j} \cdot N'}$, let $\Gamma_{j} = \textrm{Aut}(V_{j}, Q_{j})(\mathbb{Z})$ and $D_{j}$ be the period domain of polarized Hodge structures on $V_{j}$; here $(V, Q)$ is the polarized lattice associated to the period domain $D$, and $Q_{j}$ is the polarization on $V_{j}$ induced by $Q$. Then there are algebraic varieties $T_{j} \subset \Gamma_{j} \backslash D_{j}$ and algebraic maps $p_{j} : T \to T_{j}$ such that $\varphi_{j} : S \to \Gamma_{j} \backslash D_{j}$ factors through a dominant algebraic map $S \to T_{j}$, we have $\varphi_{j} = p_{j} \circ \varphi$, and the map $T \to T_{1} \times T_{2}$ induced by $p_{1}$ and $p_{2}$ is quasi-finite;
\item[(ii)] The composition $S \xrightarrow{p} T \xrightarrow{p_{1}} T_{1}$ is a Hermitian symmetric foliation of $S$. Moreover, the irreducible components occurring in the fibres of $p_{1} \circ p$ are exactly the components of the projections to $S$ of the orbits of $J_{1}$.
\end{itemize}
\end{prop}

\begin{proof}
We note the fact that $J_{j} \cdot N'$ is normal in $M$ implies that $M$ acts on $V_{j}$. The $\mathbb{Q}$-algebraic groups $J_{1} \cdot N'$ and $J_{2} \cdot N'$ are determined by the tensor invariants that they fix, so we can choose the $a_{ij}$ and $b_{ij}$ such that the representation $M \to \textrm{Aut}(V_{j}, Q_{j})$ on $V_{j}$ has kernel exactly $J_{j} \cdot N'$. The map which sends a Hodge structure on $V$ lying inside $D_{T}$ to the induced Hodge structure on $V_{j}$ induces a holomorphic map $D_{T} \to D_{j}$ which descends to a holomorphic map $r_{j} : \pi(D_{T}) \to \Gamma_{j} \backslash D_{j}$, where $\pi : D \to \Gamma \backslash D$ is the projection.\footnote{We note that $\pi(D_{T})$ is an analytic subvariety of $\Gamma \backslash D$. This can be seen by exhibiting it as the image of the map $q : (\Gamma \cap N(\mathbb{R})) \backslash D_{T} \to \Gamma \backslash D$ and using the fact that $q$ is a holomorphic map definable in the o-minimal structure $\mathbb{R}_{\textrm{alg}}$, see \cite{defpermap}.} The restrictions of the $r_{j}$ to $T$ induce Griffiths-transverse holomorphic maps $q_{j} : T \to \Gamma_{j} \backslash D_{j}$. By \cite{OMINGAGA} the maps $q_{j}$ factor as $q_{j} = \iota_{j} \circ \an{p}_{j}$ where $p_{j} : T \to T_{j}$ is a dominant map of complex algebraic varieties and $\iota_{j} : \an{T}_{j} \hookrightarrow \Gamma_{j} \backslash D_{j}$ is an embedding.

To show that the algebraic map $T \to T_{1} \times T_{2}$ is quasi-finite, it suffices to show same statement for the fibres of $q : \pi(D_{T}) \to \Gamma_{1} \backslash D_{1} \times \Gamma_{2} \backslash D_{2}$. Using the fact that $q$ is a definable morphism (see \cite{defpermap}) in the definable structure inherited from $D_{T}, D_{1}$ and $D_{2}$ this reduces to showing the fibres of $\widetilde{q} : D_{T} \to D_{1} \times D_{2}$ are finite. This is easily checked by identifying $D_{T}$ with an $N(\mathbb{R})^{\circ}$-orbit of conjugacy classes of the Deligne torus, and likewise for the image of $\widetilde{q}$, and using the fact that the map $N_{\mathbb{R}} \to N_{\mathbb{R}}/J_{1,\mathbb{R}} \times N_{\mathbb{R}} / J_{2,\mathbb{R}}$ is quasi-finite.

We now attend to part (ii). We observe that the components of the fibres of $p_{1} \circ p$ are exactly  the maximal subvarieties $Y \subset S$ on which the variation $\mathbb{V}_{1}$ is constant. Since the kernel of the representation $M \to \textrm{Aut}(V_{1}, Q_{1})$ is $J_{1} \cdot N'$, necessarily such subvarieties have algebraic monodromy contained in $J_{1}$, and are maximal for this property. It follows that these $Y$ are the components of the projections $\pi(\widetilde{\varphi}^{-1}(\mathcal{W}_{t}))$ of \autoref{hermsym}, hence the result.
\end{proof}

\section{Interlude on Compatible Families}

This is a short section concerning the construction of infinite ``compatible sequences'' of points lying in an infinite ``compatible family'' of constructible sets. It is applied in the next section to give a jet-theoretic description of the equivalence relation (E3) (over an open subset of $S$). The essential idea is that we wish to know that if we can construct a jet $j_{r} \in J^{d}_{r} S$ satisfying certain constructible algebraic conditions for every $r$, then there is an infinite sequence $\{ j_{r} \}_{r \geq 0}$, with $j_{r+1}$ projecting to $j_{r}$, satisfying these conditions at every finite stage. We will actually work in a slightly more general setup than that of jet theory for the purposes of our application.

Our setup is as follows. We suppose we have $K$-varieties $E_{r}$ for every integer $r \geq 0$ and projections $\pi^{r'}_{r} : E_{r'} \to E_{r}$ for $r' \geq r$ such that $\pi^{r'}_{r} \circ \pi^{r''}_{r'} = \pi^{r''}_{r}$ whenever the composition makes sense. Our model is of course the case where $E_{r} = J^{d}_{r} S$ for a $K$-variety $S$.

\begin{defn}
\label{compfamdef}
Suppose that $\mathcal{T}_{r} \subset E_{r}(\mathbb{C})$ for $r \geq 0$ are sets. We say that $\{ \mathcal{T}_{r} \}_{r \geq 0}$ is a \emph{compatible family} if the projections $\pi^{r'}_{r} : E_{r'} \to E_{r}$ restrict to give maps $\mathcal{T}_{r'} \to \mathcal{T}_{r}$. 
\end{defn}

\begin{defn}
\label{compsecdef}
We say that $\{ j_{r} \}_{r \geq 0}$ with $j_{r} \in E_{r}$ form a \emph{compatible sequence} if for $r' \geq r$, we have $\pi^{r'}_{r}(j_{r'}) = j_{r}$. 
\end{defn}

\begin{lem}
\label{arbextlem}
Suppose we have a compatible family $\mathcal{T}_{r} \subset E_{r}(\mathbb{C})$, and each $\mathcal{T}_{r}$ is a constructible algebraic subset. Let $j_{r_{0}} \in \mathcal{T}_{r_{0}}$ be a jet for some fixed $r_{0}$. Then if the fibre of $\pi^{r'}_{r}$ above $j_{r_{0}}$ intersects $\mathcal{T}_{r'}$ for each $r' \geq r_{0}$, then $j_{r_{0}}$ belongs to a compatible sequence $\{ j_{r} \}_{r \geq 0}$ with $j_{r} \in \mathcal{T}_{r}$ for each $r \geq 0$.
\end{lem}

\begin{proof}
We will say that $j_{r}$ admits \emph{arbitrary extensions} with respect to the family $\mathcal{T}_{r}$ if the fibre of $\pi^{r'}_{r}$ above $j_{r}$ intersects $\mathcal{T}_{r'}$ for each $r' \geq r$; i.e., if $j_{r}$ satisfies the stated hypothesis on $j_{r_{0}}$. It suffices to show that if $j_{r}$ admits arbitrary extensions with respect to the family, then there exists a jet $j_{r+1} \in \mathcal{T}_{r+1}$, projecting onto $j_{r}$ which also admits arbitrary extensions with respect to the family; given this fact we can construct the infinite sequence recursively. 

Denote by $\mathcal{K}_{r+k}$ the intersection of $\mathcal{T}_{r+k}$ with the fibre above $j_{r}$ for each positive integer $k \geq 1$. Then the condition that $j_{r}$ admits arbitrary extensions tells us that the constructible sets $\textrm{im}(\mathcal{K}_{r+k} \to \mathcal{K}_{r+1})$ for $k \geq 1$ form a descending sequence of constructible sets which is non-empty at each finite stage. A countable intersection of constructible sets is non-empty if and only if every finite intersection of the sets is non-empty. It follows that we may pick $j_{r+1} \in \bigcap_{k \geq 1} \textrm{im}(\mathcal{K}_{r+k} \to \mathcal{K}_{r+1})$ which is a jet in $\mathcal{T}_{r+1}$ extending $j_{r}$ which admits arbitrary extensions.
\end{proof}

\section{Properties of Equivalence Relations}
\label{equivrelsec}

In this section, we continue with the setup of  \autoref{hermfolsec}. We recall that \autoref{hermfolsec} reduced to the case where we had a proper period map $\varphi$. We now need to show that reducing to the case of proper period maps can be done in the category of $\Qbar$-algebraically-defined polarized variations of Hodge structure. 

\begin{lem}
\label{extvar}
Let $\mathbb{V}$ be a $\Qbar$-algebraically-defined polarized variation of Hodge structure on a smooth $\Qbar$-algebraic variety $S$. Then after possibly replacing $S$ with a finite \'etale cover\footnote{We note all finite \'etale covers of $\Qbar$-varieties admit $\Qbar$-algebraic models.}, the variation $\mathbb{V}$ extends to a $\Qbar$-algebraically-defined  polarized variation on a smooth $\Qbar$-variety $S'$ such that $S \subset S'$ is a Zariski open subset, we have $\restr{\mathbb{V}'}{S} = \mathbb{V}$, and the period map $\varphi' : S' \to \Gamma \backslash D$ is proper.
\end{lem}

\begin{proof}
By Hironaka's theorem \cite{hironaka}, we may fix a $\Qbar$-algebraic compactification $S \subset \overline{S}$ with boundary a normal crossings divisor $E$. After possibly passing to a finite \'etale cover, we may assume that $\mathbb{V}$ has unipotent local monodromy relative to $E$. Using Deligne's canonical extension \cite{canonext} $\overline{\mathcal{H}}$ we obtain a $\Qbar$-algebraic\footnote{See \cite{266936}.} extension of $\mathcal{H}$ to $\overline{S}$. There is also an extension $\overline{F}^{\bullet}$ (see \cite{filtext}) of the filtration $F^{\bullet}$ to $\overline{\mathcal{H}}$, which can then be seen to be $\Qbar$-algebraic as it coincides with the $\Qbar$-algebraic filtration $F^{\bullet}$ on a dense open $\Qbar$-algebraic subset. The maximal subvariety $S' \subset \overline{S}$ to which $\mathbb{V}$ extends to a variation of Hodge structure is then simply the maximal variety to which the complex local system $\mathbb{V}_{\mathbb{C}}$ extends, which can be seen to be the domain of definition for the maximal extension $\nabla'$ of the connection $\nabla$ by the Riemann-Hilbert correspondence. The variety $S'$ is seen to be $\Qbar$-algebraic from the fact that $\nabla$ is, and one also easily deduces that $\nabla'$ is $\Qbar$-algebraic as it agrees with $\nabla$ on $S$. That the period map $\varphi' : S' \to \Gamma \backslash D$ is proper is argued as in \cite[Cor.~13.7.6]{CMS}. Defining $\mathcal{H}' = \restr{\mathcal{H}}{S'}$ and $F'^{\bullet} = \restr{\overline{F}^{\bullet}}{S'}$, it then remains to deal with the algebraic incarnation of the polarization $\mathcal{Q}' : \mathcal{H}' \otimes \mathcal{H}' \to \mathcal{O}_{S'}$. The existence of $\mathcal{Q}'$ on the complex algebraic level can be established by first extending $\mathcal{Q}$ analytically, and using the fact that the closure of the graph of $\mathcal{Q}$ will be algebraic. One then descends $\mathcal{Q}'$ to $\Qbar$ by comparison with the $\Qbar$-algebraic polarization $\mathcal{Q} : \mathcal{H} \otimes \mathcal{H} \to \mathcal{O}_{S}$.
\end{proof}

\medskip

Using the above Lemma, we now assume that our proper period map $\varphi : S \to T \subset \Gamma \backslash D$ arises from a $K$-algebraically-defined variation of Hodge structure, where $K \subset \Qbar \subset \mathbb{C}$. We denote by $S \xrightarrow{q} U \xrightarrow{r} T$ the Stein factorization of the map $p : S \to T$; we note that $q$ and $U$ are $K$-algebraic, $q$ is proper, $U$ is normal, and $r$ is finite. We write $R_{2}, R_{3} \subset S \times S$ for the graphs of the equivalence relations (E2) and (E3).

\begin{lem}
\label{equivdesc}
For $i = 2, 3$, the image $q(R_{i})$ of $R_{i}$ under the map $S \times S \to U \times U$ is an equivalence relation defined exactly as in \autoref{equivdef} but with respect to the period map $U \to \Gamma \backslash D$ through which $\varphi$ factors.
\end{lem}

\begin{proof}
What needs to be checked is that if $u = q(s)$ and $\psi : B \to D$ is a distinguished local lift around $s \in S$, then $\psi$ factors through an open neighbourhood of $U$. This is easily checked using the fact that $U$ is normal, which implies its analytic germs are irreducible.
\end{proof}

\begin{prop}
\label{E2isnormalization}
Consider the factorization $S \xrightarrow{p^{\nu}} T^{\nu} \to T$ of $\varphi$. Then the fibres of $p^{\nu}$ are the equivalence classes of the relation (E2) described in \autoref{equivdef}.
\end{prop}

\begin{proof}
We show that the equivalence relation given by the fibres of $p^{\nu}$ is a refinement of the equivalence relation (E2), the reverse argument being analogous. Fix a point $t^{\nu} \in T^{\nu}$, and two points $s_{1}, s_{2} \in p^{\nu,-1}(t^{\nu})$. If $n : T^{\nu} \to T$ the normalization map, there is a natural bijection between the points in the fibre $n^{-1}(t)$ and the local analytic branches of $T$ at $t$; more formally we have a bijection
\[ n^{-1}(t) \xrightarrow{\sim} \{ \textrm{minimal primes in } \mathcal{O}_{\an{T}, t} \} \]
given by sending $t^{\nu} \in n^{-1}(t)$ to the kernel of the map $\mathcal{O}_{\an{T}, t} \to \mathcal{O}_{T^{\nu, \textrm{an}}, t^{\nu}}$ induced by $n$. Let us denote by $\mathcal{C}(s)$ for $s \in S$ the unique analytic component of the germ $(T, \varphi(s))$ which contains the image of the germ $(S, s)$ under $\varphi$. Then the equivalence relation given by the fibres of $p^{\nu}$ is the relation
\[ s_{1} \sim s_{2} \iff \mathcal{C}(s_{1}) = \mathcal{C}(s_{2}) . \]
If we choose distinguished local lifts $\psi_{i} : B_{i} \to D$ such that $s_{i} \in B_{i}$, the fact that $\restr{\varphi}{B_{i}} = \pi \circ \psi_{i}$ where $\pi : D \to \Gamma \backslash D$ is the natural projection then implies that there exists $\gamma \in \Gamma$ such that $\mathcal{C}(\psi_{1}, s_{1}) = \gamma \cdot \mathcal{C}(\psi_{2}, s_{2})$, and so if $s_{1}$ and $s_{2}$ lie in the same fibre of $p^{\nu}$ then $s_{1} \sim_{(E2)} s_{2}$. 
\end{proof}

\begin{lem}
\label{equidim}
The equivalence relation $q(R_{2}) \subset U \times U$ is an equidimensional subvariety.
\end{lem}

\begin{proof}
Let $r^{\nu} : U \to T^{\nu}$ be factorization of $r$ through the normalization $T^{\nu} \to T$. As $r$ is quasi-finite, so is $r^{\nu}$. The subvariety $q(R_{2})$ is isomorphic to $U \times_{T^{\nu}} U$, so we are reduced the following claim: given any quasi-finite surjective map $f : X \to Y$ of normal irreducible varieties, $X \times_{Y} X$ is equidimensional. By \cite[IV, Part 3, 14.4.4]{EGA} the map $f$ is necessarily universally open, hence so is the base change $X \times_{Y} X \to X$. Then every component of $X \times_{Y} X \to X$ dominates $X$, so the result follows.
\end{proof}

\begin{lem}
\label{closenanal}
The locus $R_{3} \subset S \times S$ is a countable union $\bigcup_{j = 1}^{\infty} C_{j}$ of closed analytic subvarieties. Moreover, choosing one of the $C_{j}$ in this union, we may describe $C_{j}$ locally as follows. If $(s_{1}, s_{2}) \in C_{j}$, there is a neighbourhood $B_{1} \times B_{2}$ of $(s_{1}, s_{2})$ such that $C_{j} \cap (B_{1} \times B_{2})$ is of the following form: there are local period maps $\psi'_{i} : B_{i} \to \widetilde{T}$ such that
\[ C_{j} \cap (B_{1} \times B_{2}) = \{ (s'_{1}, s'_{2}) : (p_{2} \circ \psi'_{1})(s'_{1}) = (p_{2} \circ \psi'_{2})(s'_{2}) \} ,\]
where $D_{T} \simeq D_{1} \times D_{2}$ is the factorization of $D_{T}$ induced by the almost-direct decomposition $N = J_{1} \cdot J_{2}$ of \autoref{hermfolprop}, and $p_{2} : D_{T} \to D_{2}$ is the projection. Furthermore, the images of the $\psi'_{i}$ may be taken to lie in an analytic set $E \times F$, where $E \subset D_{1}$ is open, and $E \times F$ is a component of $\widetilde{T} \cap V$ for some open neighbourhood $V$.
\end{lem}

\begin{proof}
It suffices to work locally on $S \times S$, so suppose that $B_{1} \times B_{2} \subset S \times S$ is a neighbourhood containing a point $(s_{1}, s_{2}) \in R_{3}$. Recall that \autoref{E3redef} tells us that there exists a countable sequence $\{ g_{j} \}_{j = 1}^{\infty}$ of elements of $G(\mathbb{C})$ such that the condition (E3) on points $(s'_{1}, s'_{2})$ is equivalent to 
\begin{equation}
\label{firstrephrase}
\textrm{ there exists }\psi'_{1}, \psi'_{2} \textrm{ and } g' \in \bigcup_{j = 1}^{\infty} H \cdot g_{j} \textrm{ such that } g' \cdot \mathcal{C}(\psi'_{1}, s'_{1}) = \mathcal{C}(\psi'_{2}, s'_{2}) .
\end{equation}
We note that, because we can always assume the elements of $\Gamma$ lie in the sequence $\{ g_{j} \}_{j = 1}^{\infty}$, we may assume that the condition on $(s'_{1}, s'_{2})$ of \autoref{firstrephrase} holds for any distinguished local lifts $\psi'_{i} : B_{i} \to D$ with $s'_{i} \in B_{i}$ if and only if it holds for a single such pair $\psi_{1}, \psi_{2}$ of lifts. It therefore suffices to consider the condition $C_{j}$ given by
\begin{equation}
\label{sndtolast}
\textrm{ there exists } g' \in H \cdot g_{j} \textrm{ such that } g' \cdot \mathcal{C}(\psi_{1}, s'_{1}) = \mathcal{C}(\psi_{2}, s'_{2}) ,
\end{equation}
for some fixed pair of distinguished lifts $\psi_{1}$ and $\psi_{2}$, and show that this condition gives an analytic condition on $B_{1} \times B_{2}$ of the prescribed form. We therefore assume without loss of generality that $(s_{1}, s_{2}) \in C_{j}$. Shrinking the $B_{i}$ if necessary, we may also assume that the smallest analytic set\footnote{That is, a closed analytic subset of an open set.} $E_{i}$ containing $\psi(B_{i})$ is locally irreducible; in particular, that for $i = 1, 2$ we have $\mathcal{C}(\psi_{i}, s'_{i}) = (E_{i}, \psi_{i}(s'_{i}))$ for every point $s'_{i} \in B_{i}$. (Here we are choosing the $\psi_{i}$ so that they are distinguished local lifts for every point $s'_{i} \in B_{i}$ simultaneously, and so that $E_{i}$ is independent of $s'_{i}$.)

\hspace{1em} We recall that by \autoref{equivQgroup} the group $H$ has the same orbits as a $\mathbb{Q}$-normal subgroup $J \subset N$ of the algebraic monodromy group $N$ of the variation $\mathbb{V}$. As in \cite[III.A]{GGK} we thus obtain a factorization $D_{T} = D_{1} \times D_{2}$ where the $J(\mathbb{R})^{\circ}$-orbits are those of the form $D_{1} \times \{ t_{2} \}$ for $t_{2} \in D_{2}$. The fact that $H$ sends germs of $\widetilde{T}$ to germs of $\widetilde{T}$ (i.e., that $\widetilde{T}$ is locally foliated by $H$-orbits) necessarily implies that $\widetilde{T} = D_{1} \times \widetilde{T}_{2}$ for some closed analytic subset $\widetilde{T}_{2} \subset D_{2}$. Consequently, any germ $(\widetilde{T}, t)$ of $\widetilde{T}$ decomposes as a product $(U, t_{1}) \times (F, t_{2})$, where $(U, t_{1}) = (D_{1}, t_{1})$. Thus, shrinking the $E_{i}$ and $B_{i}$ for $i = 1, 2$ if necessary, we may assume that $E_{i} = U_{i} \times F_{i}$, where $U_{i}$ is irreducible and open in $D_{1}$. 

\hspace{1em} As we have a point $(s_{1}, s_{2}) \in R_{3}$, it follows that there exists $g' = h \cdot g_{j}$ such that $g' \cdot (E_{1}, \psi_{1}(s_{1})) = (E_{2}, \psi_{2}(s_{2}))$. As $h$ acts trivially on the second factor, this gives in particular that $(g_{j} \cdot F_{1}, p_{2}(g_{j} \cdot \psi_{1}(s_{1}))) = (F_{2}, p_{2}(\psi_{2}(s_{2})))$. Shrinking $B_{1}$ and $B_{2}$ if necessary we may assume that $g_{j} \cdot F_{1} = F_{2}$. We also let $\psi'_{1} = g_{j} \cdot \psi_{1}$. Because $J$ (and hence $H$) acts transitively on the germs of $D_{1}$, the condition on $(s'_{1}, s'_{2}) \in B_{1} \times B_{2}$ defining $C_{j}$ becomes simply $(p_{2} \circ \psi'_{1})(s'_{1}) = (p_{2} \circ \psi_{2})(s'_{2})$, with $p_{2} : D_{T} \to D_{2}$ the projection. 

To complete the proof, take the maps $\psi'_{1}$ and $\psi_{2}$ produced in the above argument, define $\psi'_{2} = h_{0} \cdot \psi_{1}$ with $h_{0} \in H$ so that $\psi'_{1}(s_{1}) = \psi'_{2}(s_{2})$, and possibly shrinking the $B_{i}$ take $U_{1} = h_{0} \cdot U_{2}$. Note that the analytic set $E \times F$ in the statement of the proposition is $U_{1} \times F_{1}$.
\end{proof}

We will now aim to show that $R_{3}$ is in fact a $K$-algebraic subvariety of $S \times S$. For this we introduce a fourth equivalence relation $R_{4} \subset S \times S$, as follows:

\begin{defn}
We say
\begin{itemize}
\item[(E4)] that $s_{1} \sim_{(E4)} s_{2}$ if for every choice of $r \geq 0$ and $d = \dim S$ we have
\[ \eta^{d}_{r}((J^{d}_{r} S)_{s_{1}}) = \eta^{d}_{r}((J^{d}_{r} S)_{s_{2}}) , \]
where $\eta^{d}_{r}$ is as in \autoref{jetcorresp}.
\end{itemize}
\end{defn}

\begin{prop}
\label{E4prop}
The graph of (E4) satisfies the following properties:
\begin{itemize}
\item[(i)] $R_{4} \subset S \times S$ is a countable intersection of $K$-algebraic constructible sets;
\item[(ii)] we have $R_{4} \subset R_{3}$; and
\item[(iii)] there exists an open dense subset $S^{\circ} \subset S(\mathbb{C})$ such that: 
\begin{itemize}
\item[(a)] we have $R_{4} \cap (S^{\circ} \times S^{\circ}) = R_{3} \cap (S^{\circ} \times S^{\circ})$, and
\item[(b)] $R_{3} \cap (S^{\circ} \times S^{\circ})$ is dense in $R_{3}$ in the analytic topology.
\end{itemize}
\end{itemize}
\end{prop}

\begin{proof} 
~
\begin{itemize}
\item[(i)] Fix $d = \dim S$, and let $p_{r} : \mathcal{T}^{d}_{r} \to J^{d}_{r} S$ be the $G$-torsor of \autoref{torsorprop}, and $\alpha_{r} : \mathcal{T}^{d}_{r} \to J^{d}_{r} \ch{D}$ be the associated $G$-invariant map. The condition that $s_{1} \sim_{(E4)} s_{2}$ is the condition that 
\[ \alpha_{r}(p_{r}^{-1}((J^{d}_{r} S)_{s_{1}})) = \alpha_{r}(p_{r}^{-1}((J^{d}_{r} S)_{s_{2}})) \]
for all $r \geq 0$. To show (i), it suffices to show that for fixed $r$ the condition
\begin{equation}
\label{constrinc}
\alpha_{r}(p_{r}^{-1}((J^{d}_{r} S)_{s_{1}})) \subset \alpha_{r}(p_{r}^{-1}((J^{d}_{r} S)_{s_{2}})) , 
\end{equation}
gives a constructible $K$-algebraic condition on $S \times S$. To do this we first observe that the locus $\mathcal{C}_{2} \subset J^{d}_{r} \ch{D} \times S$ defined by
\[ \mathcal{C}_{2} = \{ (j, s_{2}) : j \in \alpha_{r}(p_{r}^{-1}((J^{d}_{r} S)_{s_{2}})) \} \]
is $K$-constructible, as it is the projection under $\alpha_{r} \times \textrm{id}$ of the $K$-algebraic locus $\{ (m_{r,2}, s_{2}) : m_{r,2} \in p_{r}^{-1}((J^{d}_{r} S)_{s_{2}}) \} \subset \mathcal{T}^{d}_{r} \times S$. We may then consider the locus $\mathcal{D}_{1} \subset \mathcal{T}^{d}_{r} \times S \times S$ defined by 
\[ \mathcal{D}_{1} = \{ (m_{r,1}, s_{1}, s_{2}) : p_{r}(m_{r,1}) \in (J^{d}_{r} S)_{s_{1}} \textrm{ and } (\alpha_{r}(m_{r,1}), s_{2}) \not\in \mathcal{C}_{2} \}. \]
This is once again a $K$-constructible condition. The projection of $\mathcal{D}_{1}$ to $S \times S$ defines the complement of the condition (\ref{constrinc}), which shows (i).

\item[(ii)] Suppose that $s_{1} \sim_{(E4)} s_{2}$, and let $\mathcal{A}_{r} \subset \mathcal{T}^{d}_{r} \times \mathcal{T}^{d}_{r}$ be the locus defined by
\[ \mathcal{A}_{r} = \{ (m_{r,1}, m_{r,2}) : p_{r}(m_{r,i}) \in (J^{d}_{r,nd} S)_{s_{i}} \textrm{ for }i=1,2 \textrm{ and } \alpha_{r}(m_{r,1}) = \alpha_{r}(m_{r,2}) \} . \]
The hypothesis that $s_{1} \sim_{(E4)} s_{2}$ ensures that the algebraic sets $\mathcal{A}_{r}$ form a compatible family (see \autoref{compfamdef}) with each $\mathcal{A}_{r}$ non-empty. The intersection $\bigcap_{r} \pi^{r}_{0}(\mathcal{A}_{r}) \subset \mathcal{T}^{d}_{0}$ is a countable intersection of constructible sets, non-empty at each finite stage, hence non-empty. Choosing a point $(m_{0,1}, m_{0,2})$ in this intersection, we may find by \autoref{arbextlem} a compatible sequence $\{ (m_{r,1}, m_{r,2}) \}_{r \geq 0}$ extending $(m_{0,1}, m_{0,2})$ such that $(m_{r,1}, m_{r,2}) \in \mathcal{A}_{r}$. Let $j_{r,i} = p_{r}(m_{r,i})$. Applying \autoref{torsorprop}(i) and \autoref{torsorprop}(iii) we find that $\psi_{m_{0,1}} \circ j_{r,1} = \psi_{m_{0,2}} \circ j_{r,2}$ for all $r$, where $\psi_{m_{0,i}} : B_{i} \to \ch{D}$ are the local period maps of \autoref{torsorprop}. As the jets $j_{r,i}$ are non-degenerate by construction, we are now in the situation of \autoref{sameimagelem} with $X_{i} = B_{i}$ and $Y = \ch{D}$, hence the images of the local period maps $\psi_{m_{0,1}}$ and $\psi_{m_{0,2}}$ have the same image germ (see \autoref{imdef}) at $t = \psi_{m_{0,1}}(j_{0,1}) = \psi_{m_{0,2}}(j_{0,2})$, and $s_{1} \sim_{(E3)} s_{2}$.

\item[(iii)] In we understand a map $f : X \to Y$ of analytic spaces to be \emph{smooth} if it satisfies the defining conditions of \cite[Thm 3.1, Def. 3.2]{SHC_1960-1961__13_1_A9_0}, parallel to the notion with the same name in algebraic geometry. The notion is local on the source, so makes sense at the level of germs. If $\psi = g \cdot \psi'$ is a local period map with $\psi'$ a local lift, we will also write $\mathcal{C}(\psi, s)$ for $g \cdot \mathcal{C}(\psi', s)$. Let us then consider the set
\[ S^{\circ} = \left\{ s \in S(\mathbb{C}) : (B, s) \xrightarrow{\psi} \mathcal{C}(\psi, s) \textrm{ is smooth for all local period maps } \psi \textrm{ at }s \right\} . \]
Using the usual equivalence between smoothness and formal smoothness (see \cite[Thm. 3.1]{SHC_1960-1961__13_1_A9_0}), one has (for $d = \dim S$, say)
\[ S^{\circ} = \left\{ s \in S(\mathbb{C}) : \begin{array}{c} J^{d}_{r} \psi : (J^{d}_{r} B)_{s} \to (J^{d}_{r} \, \mathcal{C}(\psi, s))_{\psi(s)} \textrm{ is surjective for all }r \\ \textrm{and any local period map } \psi \textrm{ with }s \in B  \end{array} \right\},  \]
As the action of $G(\mathbb{C})$ on germs of local period maps is transitive, it is clear that to check this condition at a point $s \in S(\mathbb{C})$ it suffices to consider the map $(B, s) \xrightarrow{\psi} \mathcal{C}(\psi, s)$ on germs for a single local period map $\psi$. Thus, we may choose local period maps $\{ \psi_{i} \}_{i \in I}$ relative to an open cover $\{ B_{i} \}_{i \in I}$, and for each $i$ let $E_{i} \subset \ch{D}$ be an analytic set such that $(E_{i}, \psi_{i}(s_{i}))$ represents the germs $\mathcal{C}(\psi_{i}, s_{i})$ for all $s_{i} \in B_{i}$. (We can find such $E_{i}$ after possibly shrinking the $B_{i}$.) Then we find that $S^{\circ} \cap B_{i}$ coincides with the smooth locus of the map $\psi_{i} : B_{i} \to E_{i}$, hence is an analytic open set.

\begin{itemize}
\item[(a)] We now suppose we have $(s_{1}, s_{2}) \in R_{3} \cap (S^{\circ} \times S^{\circ})$. Then we may find local period maps $\psi_{i} : B_{i} \to \ch{D}$ with $s_{i} \in B_{i}$ such that $\mathcal{C}(\psi_{1}, s_{1}) = \mathcal{C}(\psi_{2}, s_{2})$. As $s_{i} \in S^{\circ}$ we have
\[ \psi_{1} \circ (J^{d}_{r} S)_{s_{1}} = (J^{d}_{r} \, \mathcal{C}(\psi_{1}, s_{1}))_{\psi_{1}(s_{1})} = (J^{d}_{r} \, \mathcal{C}(\psi_{2}, s_{2}))_{\psi_{2}(s_{2})} = \psi_{2} \circ (J^{d}_{r} S)_{s_{2}} , \]
and hence $\eta^{d}_{r}((J^{d}_{r} S)_{s_{1}}) = \eta^{d}_{r}((J^{d}_{r} S)_{s_{2}})$ by the defining property of $\eta^{d}_{r}$. Hence $(s_{1}, s_{2}) \in R_{4} \cap (S^{\circ} \times S^{\circ})$, and the result follows.
\item[(b)] Fix a point $(s_{1}, s_{2}) \in R_{3}$. Applying \autoref{closenanal} we may reduce to the following setting: we have local period maps $\psi_{i} : B_{i} \to E \times F \subset \widetilde{T}$ with $(s_{1}, s_{2}) \in B_{1} \times B_{2}$, a locus $C \subset B_{1} \times B_{2}$ defined on $(s'_{1}, s'_{2}) \in B_{1} \times B_{2}$ by $(p_{2} \circ \psi_{1})(s'_{1}) = (p_{2} \circ \psi_{2})(s'_{2})$, and we need to show that $C \cap (S^{\circ} \times S^{\circ})$ is dense in $C$. By factoring the $\psi_{i}$ through the map $S \xrightarrow{q} U$ and replacing $S$ by $U$ we may assume the $\psi_{i}$ are quasi-finite. By the Open Mapping Theorem \cite[pg.107]{Grauert1984} we obtain that the maps $\psi_{i} : B_{i} \to E \times F$ are open, and so shrinking $E \times F$ and the $B_{i}$ if necessary we may assume they are surjective. The base-change $C \to \Delta_{2}$ of $\psi_{1} \times \psi_{2}$ is similarly open and surjective, where $\Delta_{2} \subset E \times F \times E \times F$ defined by 
\[ \Delta_{2} = \{ (t_{1}, t_{2}, t'_{1}, t'_{2}) : t_{2} = t'_{2} \} . \]
\hspace{1em} Since the inverse image of a dense set under an open map is dense, it suffices to find a dense open subset $T^{\circ} \subset E \times F$ such that $\psi_{i}^{-1}(T^{\circ}) \subset S^{\circ}$ for $i = 1, 2$ and $\Delta_{2} \cap (T^{\circ} \times T^{\circ})$ is dense in $\Delta_{2}$. We take  
\[ T^{\circ} = \left\{ (e, f) \in (E \times F) :  \begin{array}{c} \psi_{i} : (B_{i}, s) \to (E \times F, (e, f)) \textrm{ is smooth } \\ \textrm{for all }s \in \psi_{i}^{-1}(e, f) \textrm{ and }i = 1,2 \end{array} \right\} . \] 
As $T^{\circ}$ is the complement of an analytic subset of $E \times F$ and $E$ is irreducible, it follows that the subset $F^{\circ} \subset F$ where the fibres of $T^{\circ} \to F$ have dense image in $E$ (i.e., points $f \in F$ where  $\dim \, ((E \times F) \setminus T^{\circ}) \cap (E \times \{ f \}) < \dim E$) is a dense subset of $F$.

\hspace{1em} Now let $(t_{1}, t_{2}, t'_{1}, t_{2})$ be a point of $\Delta_{2}$. By the density of $F^{\circ}$ in $F$, we may choose a convergent sequence $t^{i}_{2} \to t_{2}$ with $t^{i}_{2} \in F^{\circ}$ for all $i$. As each fibre $T^{\circ}_{t^{i}_{2}}$ of the projection $T^{\circ} \to F$ is dense in $E \times \{t^{i}_{2}\}$, we may then pick sequences $t^{i}_{1} \to t_{1}$ and $t'^{i}_{1} \to t'_{1}$ such that $(t^{i}_{1}, t^{i}_{2}, t'^{i}_{1}, t^{i}_{2}) \in (T^{\circ} \times T^{\circ}) \cap \Delta_{2}$ for all $i$. The result follows.
\end{itemize}
\end{itemize}
\end{proof}

\begin{lem}
\label{R3alg}
$R_{3} \subset S \times S$ is a $K$-algebraic subvariety.
\end{lem}

\begin{proof}
From \autoref{closenanal} and \autoref{E4prop} we learn the following: the locus $R_{3}$ is the closure in the analytic topology of a locus $R_{4} \subset R_{3}$ which is a countable intersection of $K$-constructible sets. Any such intersection is necessarily of the form
\[ \bigcup_{i = 1}^{n} \left( V_{i} \cap \bigcap_{j = 1}^{\infty} U_{i,j} \right) , \]
where $V_{i}$ is a closed $K$-algebraic set and the $U_{i,j}$ are dense $K$-algebraic Zariski open subsets of $V_{i}$. As a countable intersection of Zariski open sets remains dense (even in the analytic topology), it follows that $R_{3} = V_{1} \cup \cdots \cup V_{n}$.
\end{proof}

\section{Main Results}
\label{mainresultssec}

\subsection{The Quasi-finite Case}

We now give a stronger version of \autoref{E3prop}(ii).

\begin{prop}
\label{finiteorherm}
Either the map
\[ S(\mathbb{C}) \big/ \sim_{(E2)} \hspace{0.5em} \to  \hspace{0.5em} S(\mathbb{C}) \big/ \sim_{(E3)} \]
induced by the identity has finite fibres, or the canonical Hermitian symmetric foliation of $T$ in \autoref{hermfolprop}(ii) is positive-dimensional (c.f. \autoref{hermfoldef}).
\end{prop}

\begin{proof}
Using the fact that (E2) is given by the fibres of $S \to T^{\nu}$ we may reduce to the same statement for $U(\mathbb{C})$. As $U \to T^{\nu}$ is quasi-finite, so are the equivalence classes under $q(R_{2})$, so the claim reduces to showing that if some equivalence class of $U(\mathbb{C})$ under (E3) is infinite, then the Hermitian symmetric foliation on $T$ given by \autoref{hermfolprop}(ii) is positive-dimensional. By \autoref{R3alg}, we learn that each equivalence class of (E3) in $U$ is an algebraic subvariety, hence if some equivalence class is not finite, it must be uncountable. In this case the set $\widetilde{H}$ of \autoref{E3redef} responsible for inducing equivalences under (E3) must also be uncountable, and hence the group $H = J_{1}$ of \autoref{hermfolsec} which induces the Hermitian symmetric foliation must be positive-dimensional.
\end{proof}

\begin{lem}
\label{imagedesc}
Suppose that $X$ is a $K$-algebraic variety, that $f_{\mathbb{C}} : X_{\mathbb{C}} \to Y$ is a surjective map of irreducible complex algebraic varieties with $Y$ normal, and the graph $R \subset X_{\mathbb{C}} \times X_{\mathbb{C}}$ of the equivalence relation induced by the fibres of $f_{\mathbb{C}}$ is stable under $\textrm{Aut}(\mathbb{C}/K)$. Then $Y$ admits a $K$-algebraic model for which $f_{\mathbb{C}}$ is the base-change of a $K$-algebraic map $f : X \to Y$.
\end{lem}

\begin{proof}
Let $\sigma \in \textrm{Aut}(\mathbb{C}/K)$. Consider the projection $C^{\sigma}$ of $R$ to $Y \times Y^{\sigma}$. We observe that $C^{\sigma}$ is a bijective correspondence from $Y$ to $Y^{\sigma}$; indeed, it amounts to the map $t \mapsto f_{\mathbb{C}}^{\sigma}(f_{\mathbb{C}}^{-1}(t))$, and the fact that $R$ is stable under $\textrm{Aut}(\mathbb{C}/K)$ means that $f_{\mathbb{C}}$ and $f^{\sigma}_{\mathbb{C}}$ have the same sets of fibres. As $Y$ (and hence $Y^{\sigma}$) is normal, the correspondences $C^{\sigma}$ induce isomorphisms $\eta^{\sigma} : Y \to Y^{\sigma}$. (For this one may use that the projections $C^{\sigma} \to Y$ and $C^{\sigma} \to Y^{\sigma}$ are bijective algebraic maps to normal varieties, hence isomorphisms by \cite[Prop. 14.7]{normalization}.) The isomorphisms $\eta^{\sigma}$ satisfy $(\eta^{\sigma})^{\tau} \circ \eta^{\tau} = \eta^{\sigma \tau}$ for $\sigma, \tau \in \textrm{Aut}(\mathbb{C}/K)$, hence form a descent datum. One therefore obtains $Y = Y'_{\mathbb{C}}$ for a $K$-variety $Y'$, and one easily checks that the graph of $f_{\mathbb{C}}$ is $\textrm{Aut}(\mathbb{C}/K)$-invariant with respect to the $K$-structure on the product $X \times Y'$.
\end{proof}

\begin{lem}
\label{trivfoilproof}
Suppose that the Hermitian symmetric foliation of \autoref{hermfolprop}(ii) is trivial.  Then $T$ and $p^{\nu} : S \to T^{\nu}$ are defined over $\overline{\mathbb{Q}}$.
\end{lem}

\begin{proof}
By \autoref{equidim} we know that $q(R_{2})$ is equidimensional, that $q(R_{2}) \subset q(R_{3})$, and by \autoref{R3alg} that $q(R_{3})$ is defined over $\Qbar$. It follows from \autoref{finiteorherm} that if the canonical Hermitian symmetric foliation of \autoref{hermfolprop}(ii) is trivial, then $q(R_{3})$ must have the same dimension as $q(R_{2})$. But then $q(R_{2})$ is a union of finitely many $\Qbar$-algebraic components of $q(R_{3})$, hence $\Qbar$-algebraic. We get the same result for $R_{2}$ using the fact that $R_{2} = q^{-1}(q(R_{2}))$. The result then follows from \autoref{imagedesc} with $R = R_{2}$, $X = S$, $Y = T^{\nu}$ and $f_{\mathbb{C}} = p^{\nu}$.
\end{proof}

\subsection{The Hermitian Symmetric Case}

In this section we work in the setting of \autoref{hermfolprop}. As \autoref{trivfoilproof} above will prove the main theorems in the case where the Hermitian symmetric foliation of \autoref{hermfolprop} is trivial, it remains to deal with the Hermitian symmetric case directly. This will essentially amount to the well-known fact that if $q : U \to \Gamma' \backslash D'$ is a dominant period map to a Shimura variety $X' = \Gamma' \backslash D'$ from a $\Qbar$-variety $U$, then $q$ admits a $\Qbar$-algebraic model. However we need to deal with the additional subtly that $T_{2} \subset \Gamma_{2} \backslash D_{2}$ is not quite a Shimura variety, but only the finite image of one under a map of Hodge varieties.

\begin{lem}
\label{T2Qbar}
There exists a Shimura variety $X' = \Gamma' \backslash D'$ and a surjective morphism $b: X' \to T_{2}$ such that:
\begin{itemize}
\item[(i)] the normalization $T^{\nu}_{2}$ of the variety $T_{2} \subset \Gamma_{2} \backslash D_{2}$ of \autoref{hermfolprop} admits a $\Qbar$-algebraic model;
\item[(ii)] the induced map $b^{\nu} : X' \to T^{\nu}_{2}$ is $\Qbar$-algebraic, where $X'$ has its canonical $\Qbar$-structure.
\end{itemize}
\end{lem}

\begin{proof}
Let us recall the construction of $T_{2}$. We denote by $M$ the Mumford-Tate group of $\mathbb{V}$, which regard as a subgroup of $\textrm{Aut}(V, Q)$ where $(V, Q)$ is a fixed polarized lattice, and let $J_{2}$ and $N'$ be as in \autoref{hermfolprop}. We recall that \autoref{hermfolprop} associated to $\mathbb{V}$ an auxiliary variation $\mathbb{V}_{2}$ and a representation $\rho : M \to M' \subset \textrm{Aut}(V_{2}, Q_{2})$, where we identify $M' = M/(J_{2} \cdot N')$, and $(V_{2}, Q_{2})$ was a polarized lattice constructed from $(V, Q)$. In particular, $T \to T_{2}$ was induced by the map $\widetilde{T} \to D_{2}$ obtained by restricting the composition
\begin{align*} 
\underbrace{\{ \textrm{reps. } \mathbb{S} \to M_{\mathbb{R}} \}\Big/\textrm{conj.}}_{D_{M}} &\xrightarrow{\alpha} \underbrace{\{ \textrm{reps. }\mathbb{S} \to M'_{\mathbb{R}} \} \Big/\textrm{conj.}}_{D_{M'}} \xrightarrow{\beta} \underbrace{\{ \textrm{reps. }\mathbb{S} \to \textrm{Aut}(V_{2}, Q_{2})_{\mathbb{R}} \} \Big/ \textrm{conj.}}_{D_{2}}
\end{align*}
where $\mathbb{S}$ is the Deligne torus.

Let us observe that the restriction of $\alpha$ to $\widetilde{T}$ is surjective onto a component $\widetilde{X}' \subset D_{M'}$; indeed, the components of $D_{M'}$ are generated by orbits of $J_{1}(\mathbb{R})^{\circ}$, the map $\alpha$ is $J_{1}(\mathbb{R})^{\circ}$-invariant, and $\widetilde{T}$ contains such orbits using the fact that orbits of $J_{1}$ agree with orbits of $H$. It follows that after choosing a level subgroup $\Gamma' \subset M'(\mathbb{Q}) \cap \textrm{Aut}(V_{2}, Q_{2})(\mathbb{Z})$, the variety $T_{2}$ is a component of the image of the map $X' = \Gamma' \backslash D_{M'} \xrightarrow{b} \Gamma_{2} \backslash D_{2}$. It therefore suffices to show that the normalization $X^{\nu}$ of the image $X = \textrm{im}(b)$ admits a $\Qbar$-algebraic model.

As the orbits of $J_{1}(\mathbb{R})^{\circ}$ are Hermitian symmetric domains and the map $D_{M} \to D_{M'}$ preserves Griffiths transverse directions, it follows that $(M', \widetilde{X}')$ is a Shimura datum, and $X'$ is a Shimura variety. \cite{OMINGAGA} then shows that $X$ and the map $b$ admit complex algebraic models. Moreover, $X'$ has a canonical $\Qbar$-algebraic model characterized by the property that the dense (in both the Zariski and analytic toplogies) set of CM points is stable under the $\textrm{Aut}(\mathbb{C}/\Qbar)$ action. If we then denote by $R'_{2} \subset X' \times X'$ the graph induced by the map $X' \to X^{\nu}$, we see that $R'_{2}$ is stable under the $\textrm{Aut}(\mathbb{C}/\Qbar)$-action. The result then follows from \autoref{imagedesc}.
\end{proof}

\begin{lem}
\label{r2Qbar}
The natural factorization $r^{\nu}_{2} : S \to T^{\nu}_{2}$ of the map $r_{2} : S \xrightarrow{p} T \xrightarrow{p_{2}} T_{2}$ appearing in \autoref{hermfolprop} admits a $\Qbar$-algebraic model.
\end{lem}

\begin{proof}
By \autoref{T2Qbar} we may assume that both $S$ and $T_{2}^{\nu}$ admit $\Qbar$-algebraic models. The result will now follow from a spreading out argument (c.f. \cite{peters}) if we can show that the map $r^{\nu}_{2}$ admits no non-trivial deformations (fixing the source and target). As any non-trivial deformation of $r^{\nu}_{2}$ induces a non-trivial deformation of $r_{2}$, it suffices to show this for $r_{2}$. The rigidity of $r_{2}$ may be checked after a base-change, so we may reduce to the Shimura case using \autoref{T2Qbar}(ii), where the result follows by the main result of \cite{dominanthyperbolic} (c.f. \cite[Ex. 2]{dominanthyperbolic}).
\end{proof}

\subsection{The General Case}

\begin{thm}
\label{bigthm}
Let $\mathbb{V}$ be any $\overline{\mathbb{Q}}$-algebraically-defined variation of Hodge structure on a smooth $\overline{\mathbb{Q}}$-algebraic variety $S$. Let $\varphi : S \to \Gamma \backslash D$ be the period map, with $D$ the period domain of Hodge structures on an integral lattice $V$ polarized by $Q$, and $\Gamma = \textrm{Aut}(V, Q)(\mathbb{Z})$. Consider the factorization $\varphi = \iota \circ \an{p}$ of the main theorem in \cite{OMINGAGA}, where $p : S \to T$ is a dominant map of algebraic varieties, and $\iota : \an{T} \hookrightarrow \Gamma \backslash D$ is a closed embedding. Then the normalization $T^{\nu}$ of $T$ admits the structure of a $\overline{\mathbb{Q}}$-algebraic variety such that the map $S \to T^{\nu}$ induced by $p$ is a $\overline{\mathbb{Q}}$-algebraic map.
\end{thm}

\begin{proof}
We argue by induction on $\dim T$. Applying \autoref{extvar} to extend the variation on $S$ we may assume that $\varphi$ is proper and so that $p$ is surjective. When the canonical Hermitian symmetric foliation of \autoref{hermfolprop} is trivial we are done by \autoref{trivfoilproof}. Otherwise \autoref{hermfolprop} gives us a non-trivial factorization $S \to T \to T_{1} \times T_{2}$. We note that both the variations $\mathbb{V}_{i}$ associated to the periods maps $r_{i} : S \to T_{i}$ admit $\Qbar$-algebraic models: indeed, if $\mathbb{V}$ admits a $\Qbar$-algebraic model then so do tensor powers, direct sums, duals and subvariations (for subvariations this is shown as in \cite[Lemma 3.3]{fieldsofdef}). We thus have that the map $r_{1}^{\nu} : S \to T^{\nu}_{1}$ admits a $\Qbar$-algebraic model by the induction hypothesis, and the map $r_{2}^{\nu} : S \to T^{\nu}_{2}$ admits a $\Qbar$-algebraic model by \autoref{r2Qbar}. It follows that the equivalence relation $R'_{2}$ defined by the fibres of $S \to T^{\nu}_{1} \times T^{\nu}_{2}$ is $\Qbar$-algebraic. Using the Stein factorization $S \xrightarrow{q} U \xrightarrow{r} T$ of $p$ we obtain the same fact for the equivalence relation $q(R'_{2}) \subset U \times U$ obtained as the image under $q \times q$. The quasi-finiteness of the composition $U \to T \to T_{1} \times T_{2}$ then implies as in \autoref{equidim} that $q(R_{2}) \subset q(R'_{2})$ is a union of components of $q(R'_{2})$, hence $\Qbar$-algebraic, and using that $R_{2} = q^{-1}(q(R_{2}))$ we learn the same fact for $R_{2}$. Arguing as in \autoref{trivfoilproof} the result follows by \autoref{imagedesc}.
\end{proof}

\begin{cor}
In the setting of \autoref{bigthm}, suppose that $\varphi : S \to T \subset \Gamma \backslash D$ is quasi-finite. Let $E \subset S$ be the inverse image of the non-unibranch locus of $T$. Then if $s \in S(\Qbar) \setminus (E(\mathbb{C}) \cap S(\Qbar))$ there does not exist $s' \in S(\mathbb{C}) \setminus S(\Qbar)$ such that $\mathbb{V}_{s} \simeq \mathbb{V}_{s'}$ as polarized integral Hodge structures.
\end{cor}

\begin{proof}
By \autoref{E2isnormalization}, the statement reduces to showing that $s$ and $s'$ cannot lie in the same fibre of the map $S \to T^{\nu}$. This follows immediately from the result of \autoref{bigthm} that $S \to T^{\nu}$ admits a $\Qbar$-algebraic model compatible with the $\Qbar$-structure on $S$.
\end{proof}

Finally, let us note that although we have worked throughout the paper with the lattice $\Gamma = G(\mathbb{Z})$, our methods and results do not seriously change if one replaces $\Gamma$ with a lattice $\Gamma' \subset \Gamma$ containing the image of monodromy, provided one knows that period images in $\Gamma' \backslash D$ are algebraic. At the time of writing, \cite{OMINGAGA} claims this when the index $[\Gamma : \Gamma']$ is finite. In general, we have the following result:

\begin{cor}
\label{anylattice}
\autoref{bigthm} continues to hold if $\Gamma$ is replaced by any subgroup $\Gamma' \subset \Gamma$ containing the monodromy group $\Gamma_{S}$ provided that the map $\varphi' : S \to \Gamma' \backslash D$ factors as $\varphi' = \iota' \circ \an{p'}$, where $p' : S \to T'$ is a dominant map of algebraic varieties, and with $\iota' : \an{T'} \hookrightarrow \Gamma' \backslash D$ a closed embedding, as in the main theorem of \cite{OMINGAGA}.
\end{cor}

\begin{proof}
The map $p'^{\nu} : S \to T'^{\nu}$ induced by $\varphi'$ factors through the map $p^{\nu} : S \to T^{\nu}$, so the equivalence relation $R_{2}' \subset S \times S$ defined by the fibres of $p'^{\nu}$ lies inside the equivalence relation $R_{2} \subset S \times S$ defined by the fibres of $p^{\nu}$. Arguing as in \autoref{bigthm} we may reduce to the proper case, and taking the Stein factorization of of $p'$, to the case of quasi-finite maps $r : U \to T^{\nu}$ and $r' : U \to T'^{\nu}$. Arguing analogously to \autoref{equidim} we may assume that the components of $R_{2}'$ are among the components of $R_{2}$, so as $R_{2}$ is $\Qbar$-algebraic, the result follows by applying \autoref{imagedesc} as before.
\end{proof}

\bibliographystyle{alpha}
\bibliography{normalization}

\end{document}